\documentclass[11pt, one side, article]{memoir}

\settrims{0pt}{0pt} 
\settypeblocksize{*}{36.5pc}{*} 
\setlrmargins{*}{*}{1} 
\setulmarginsandblock{.98in}{.98in}{*} 
\setheadfoot{\onelineskip}{2\onelineskip} 
\setheaderspaces{*}{1.5\onelineskip}{*} 
\checkandfixthelayout

\usepackage{amsthm}
\usepackage{mathtools}

\usepackage[inline]{enumitem}
\usepackage{ifthen}
\usepackage[utf8]{inputenc} 
\usepackage{xcolor}

\usepackage[backend=biber, backref=true, maxbibnames = 10, style = alphabetic]{biblatex}
\usepackage[bookmarks=true, colorlinks=true, linkcolor=blue!50!black,
citecolor=orange!50!black, urlcolor=orange!50!black, pdfencoding=unicode]{hyperref}
\usepackage[capitalize]{cleveref}

\usepackage{tikz}

\usepackage{amssymb}
\usepackage{newpxtext}
\usepackage[varg,bigdelims]{newpxmath}
\usepackage{mathrsfs}
\usepackage{dutchcal}
\usepackage{mathalfa}
\usepackage{fontawesome}
\usepackage{ebproof}
\usepackage{stmaryrd}
\usepackage{ebproof}
\usepackage{graphicx}

  \crefformat{enumi}{\card#2#1#3}
  \crefalias{chapter}{section}

  \hypersetup{final}

  \setlist{nosep}
  \setlistdepth{6}

  \usetikzlibrary{ 
  	cd,
  	math,
  	decorations.markings,
		decorations.pathreplacing,
  	positioning,
  	arrows.meta,
  	shapes,
		shadows,
		shadings,
  	calc,
  	fit,
  	quotes,
  	intersections,
    circuits,
    circuits.ee.IEC
  }
  
  \tikzset{
biml/.tip={Glyph[glyph math command=triangleleft, glyph length=.95ex]},
bimr/.tip={Glyph[glyph math command=triangleright, glyph length=.95ex]},
}

\tikzset{
	tick/.style={postaction={
  	decorate,
    decoration={markings, mark=at position 0.5 with
    	{\draw[-] (0,.4ex) -- (0,-.4ex);}}}
  }
} 
\tikzset{
	slash/.style={postaction={
  	decorate,
    decoration={markings, mark=at position 0.5 with
    	{\draw[-] (.3ex,.3ex) -- (-.3ex,-.3ex);}}}
  }
}

\theoremstyle{definition}
\newtheorem{definitionx}{Definition}[chapter]

\theoremstyle{plain}

\newtheorem{warning}[definitionx]{Warning}
\newtheorem*{theorem*}{Theorem}
\newtheorem*{proposition*}{Proposition}
\newtheorem*{corollary*}{Corollary}
\newtheorem*{lemma*}{Lemma}
\newtheorem*{warning*}{Warning}

\DeclareSymbolFont{stmry}{U}{stmry}{m}{n}
\DeclareMathSymbol\fatsemi\mathop{stmry}{"23}

\DeclareFontFamily{U}{mathx}{\hyphenchar\font45}
\DeclareFontShape{U}{mathx}{m}{n}{
      <5> <6> <7> <8> <9> <10>
      <10.95> <12> <14.4> <17.28> <20.74> <24.88>
      mathx10
      }{}
\DeclareSymbolFont{mathx}{U}{mathx}{m}{n}
\DeclareFontSubstitution{U}{mathx}{m}{n}
\DeclareMathAccent{\widecheck}{0}{mathx}{"71}

\ExplSyntaxOn
\NewDocumentEnvironment{sequation}{O{\fontsize{15pt}{15pt}\selectfont
}b}
 {
  \yufip_sequation:nnn {equation}{#1}{#2}
 }{}
\NewDocumentEnvironment{sequation*}{O{\fontsize{16pt}{16pt}\selectfont
}b}
 {
  \yufip_sequation:nnn {equation*}{#1}{#2}
 }{}
\cs_new_protected:Nn \yufip_sequation:nnn
 {
  \begin{#1}
  \mbox{#2$\displaystyle#3$}
  \end{#1}
 }
\ExplSyntaxOff

\renewcommand{\ss}{\subseteq}

\DeclarePairedDelimiter{\present}{\langle}{\rangle}

\DeclareMathOperator{\Hom}{Hom}

\DeclareMathOperator*{\colim}{colim}

\DeclareMathOperator{\ob}{Ob}

\newcommand{\iHom}{\ul{\Hom}}

\newcommand{\cat}[1]{\mathcal{#1}}
\newcommand{\Cat}[1]{\mathbf{#1}}
\newcommand{\fun}[1]{\mathrm{#1}}
\newcommand{\Fun}[1]{\mathsf{#1}}

\newcommand{\id}{\mathrm{id}}
\newcommand{\then}{\mathbin{\fatsemi}}

\newcommand{\too}{\longrightarrow}
\newcommand{\tto}{\rightrightarrows}
\newcommand{\To}[2][]{\xrightarrow[#1]{\tn{$#2$}}}

\newcommand{\from}{\leftarrow}

\newcommand{\From}[1]{\xleftarrow{#1}}

\newcommand{\card}{\,^{\#}}

\newcommand{\op}{^\tn{op}}

\newcommand{\tn}[1]{\textnormal{#1}}
\newcommand{\ol}[1]{\overline{#1}}
\newcommand{\ul}[1]{\underline{#1}}

\newcommand{\lin}[1]{\hspace{1pt}\ol{\hspace{-1pt}#1\hspace{-1pt}}\hspace{1pt}}

\newcommand{\nn}{\mathbb{N}}

\newcommand{\smset}{\Cat{Set}}
\newcommand{\grp}{\Cat{Grp}}

\newcommand{\catsharp}{\Cat{Cat}^{\sharp}}

\newcommand{\sspan}{\mathbb{S}\Cat{pan}}
\newcommand{\End}{\Cat{End}}
\newcommand{\Aut}{\Cat{Aut}}

\newcommand{\List}{\Fun{list}}
\newcommand{\lott}{\Fun{lott}}

\newcommand{\yon}{{\mathcal{y}}}
\newcommand{\poly}{\Cat{Poly}}
\newcommand{\Span}{\Cat{Span}}

\newcommand{\cart}{\tn{cart}}
\newcommand{\tocart}{\To{\cart}}
\newcommand{\polycart}{\poly^\cart}

\newcommand{\0}{\textsf{0}}

\newcommand{\tri}{\mathbin{\triangleleft}}
\newcommand{\irt}{\mathbin{\triangleright}}
\newcommand{\at}[2]{\left[#1\mathbin{/}#2\right]}

\newcommand{\indep}{\Fun{Indep}}

\newcommand{\cofree}{\mathfrak{c}}
\newcommand{\free}{\mathfrak{m}}
\newcommand{\freecol}{\mathfrak{m}^\otimes}
\newcommand{\uu}{\List}

\newcommand{\biglens}[2]{
     \begin{bmatrix}{\vphantom{f_f^f}#2} \\ {\vphantom{f_f^f}#1} \end{bmatrix}
}
\newcommand{\littlelens}[2]{
     \begin{bsmallmatrix}{\vphantom{f}#2} \\ {\vphantom{f}#1} \end{bsmallmatrix}
}
\newcommand{\lens}[2]{
  \relax\if@display
     \biglens{#1}{#2}
  \else
     \littlelens{#1}{#2}
  \fi
}

\newcommand{\indexcoclscale}[1]{\scalebox{.7}{#1}}
\newcommand{\cocl}[1]{
	\scriptsize\overset{\,\indexcoclscale{$#1$}}{\frown}\normalsize
}
\newcommand{\hyper}[1]{
	\begin{tikzpicture}[y=.5cm, font=\scriptsize, baseline=(base)]
		\node[rotate=-15] (ar) {$\nearrow$};
		\coordinate[below=3pt] (base) at (ar);
		\node[above right=-2pt and 1pt of ar.west] (f) {\indexcoclscale{$#1$}};
	\end{tikzpicture}
}

\newcommand{\upup}{\uparrow\!\uparrow}

\newcommand{\hh}[2][]{#1 \tn{#2} #1}
\newcommand{\qqand}{\hh[\qquad]{and}}

\newcommand{\hi}[4][]{#1 #2 \tn{\textit{#4}} #3}
\newcommand{\where}[1][,]{\hi[#1]{\qquad}{\quad}{where}}

\newcommand{\qqor}{\hh[\qquad]{or}}

\newcommand{\OR}{\curlyvee}
\newcommand{\garner}[1]{\mathbin{\ul{#1}}}

\newcommand{\coalg}{\tn{-}\Cat{Coalg}}
\newcommand{\ext}{\fun{Ext}}

\newcommand{\thanksAFOSR}[1]{This material is based upon work supported by the Air Force Office of Scientific Research under award numbers #1}

\allowdisplaybreaks
\setsecnumdepth{section}
\settocdepth{section}
\begin{document}

\title{A summary of categorical structures in $\poly$}

\author{David I. Spivak}


\maketitle

\begin{abstract}
In this document, we collect a list of categorical structures on the category $\poly$ of polynomial functors. There is no implied claim that this list is in any way complete. It includes: infinitely-many symmetric monoidal closed structures that distribute over $+$, an asymmetric monoidal product, various duoidalities and coclosures, as well as a discussion of monoids and comonoids for some of these. We also discuss various adjunctions of which $\poly$ is a part, including the free monad and cofree comonad and their interaction with various monoidal structures. 
\end{abstract}

\tableofcontents*
\bigskip

This document is only meant as a handy guide to the abundance of structure in $\poly$. In particular, we do not supply proofs, though we have written about most of these structures elsewhere; see \cite{spivak2021functorial}, \cite{niu2022poly}, and \cite{spivak2022polynomial}. For everything else written here, one can consider it to be only conjecture, since in some cases we have not checked all the details. Hence, if someone proves something written here---something which has not been proven elsewhere---that person should be taken to have the professional ``priority'' and credit. In particular, we wish to claim no credit for originality of anything contained in this document, though most of it was discovered independently by the author, so we also do not supply additional references. We also make absolutely no claim of completeness.

\chapter{Background and notation}

A polynomial functor $p\colon\smset\to\smset$ is any functor that's isomorphic to a coproduct of representables
\[
p\coloneqq\sum_{I: p(1)}\yon^{p[I]}.
\]
We denote by $\poly$ the category of polynomial functors and all natural transformations between them.

We will typically use the above notation---which we call \emph{standard form}---appropriately modified for $p'$, $q$, etc., e.g.
\[
p'\coloneqq\sum_{I': p'(1)}\yon^{p'[I']}
\qqor
q\coloneqq\sum_{J: q(1)}\yon^{q[J]}.
\]
Putting polynomials in this form lets us consider them as combinatorial objects. We refer to elements of $p(1)$ as \emph{positions} of $p$ and, for each $I: p(1)$, we refer to the elements of $p[I]$ as \emph{directions} at $I$.

The functor $\smset\to\smset$ corresponding to $p:\poly$ is called its \emph{extension}, and denoted
\begin{align}
\ext\colon\poly&\to\End(\smset)\label{ext_functor}\\
p&\mapsto X\mapsto\sum_{I:p(1)}X^{p[I]}
\end{align} 
The functor $\ext$ preserves coproducts, limits, and cartesian closure, and it is fully faithful, meaning that a morphism between polynomials is a natural transformation $\varphi\colon \ext(p)\to \ext(q)$. What we denote $p(1)$ is actually $\ext(p)(1)$. By the Yoneda lemma and universal property of coproducts, a morphism in $\poly$ consists of a function $\varphi_1\colon p(1)\to q(1)$ and, for each $I: p(1)$ a function $\varphi^\sharp_I\colon q[\varphi_1(I)]\to p[I]$. A map $\varphi$ is called \emph{vertical} if $\varphi_1$ is identity on positions; it is called \emph{cartesian} if for each $I: p(1)$ the function $\varphi_I^\sharp$ is a bijection; see \cref{chap.bifib} for more on this. 

The subcategory of polynomial functors and cartesian morphisms between them is denoted $\polycart$. The wide, faithful functor $\polycart\to\poly$ creates all limits, coproducts, and filtered colimits; $\polycart$ does not have coequalizers.%
\footnote{For example, $\polycart$ does not have a coequalizer of the two distinct maps $\yon^2\to\yon^2$.} 
In particular, the composite functor $\polycart\to\poly\To{\ext}\End(\smset)$ preserves all limits, coproducts, and filtered colimits. The pushout of a cartesian map along any map is again cartesian, and $\ext$ preserves such pushouts.

\chapter{Coproducts and distributive monoidal structures}

The category $\poly$ has coproducts, given by the following formula:
\begin{equation}
p+q\coloneqq\sum_{I: p(1)}\yon^{p[I]}+\sum_{J: q(1)}\yon^{q[J]}
\end{equation}
If one wants the formula to be in standard form, use case logic in the exponent:
\begin{equation}
p+q\cong\sum_{X: p(1)+q(1)}\yon^{\fontsize{8pt}{8pt}\selectfont
	\begin{cases}
  	p[X]&\tn{ if }X\in p(1)\\
  	q[X]&\tn{ if }X\in q(1)
	\end{cases}\normalsize
	}
\end{equation}

The category $\poly$ is infinitary-extensive: given a map $p\to\sum_{i:I}q_i$ for any $p:\poly$, $I:\smset$, and $q\colon I\to\poly$, there is an isomorphism
\begin{equation}
	\sum_{i:I}p_i\To{\cong}p,
\end{equation} 
where $p_i\coloneqq p\times_{\sum_{i:I}q_i}q_i$ is the pullback.

For any symmetric monoidal product $(\mathbb{I},\cdot)$ on $\smset$, there is a corresponding symmetric monoidal structure $(\yon^\mathbb{I},\odot)$ on $\poly$, where the monoidal product given as follows:
\begin{equation}\label{eqn.day_all}
p\odot q\coloneqq\sum_{(I,J): p(1)\times q(1)}\yon^{p[I]\cdot q[J]}.
\end{equation}
It always distributes over $+$:
\begin{equation}\label{eqn.day_distributes}
p\odot(q_1+q_2)\cong (p\odot q_1)+(p\odot q_2).
\end{equation}
The symmetric monoidal structure $\odot$ on $\poly$ is the Day convolution of the $\cdot$ structure on $\smset$. This means we have maps
\begin{equation}\label{eqn.Day}
	p(A)\times q(B)\to (p\odot q)(A\cdot B)
\end{equation}
natural in $A,B:\smset$. 

For any set $S$, there is a monoidal structure $(0,\vee_S)$ on $\smset$,\footnote{I learned the $\vee_1$ monoidal structure (pronounced ``or'') on $\smset$ from Richard Garner, and Solomon Bothwell later informed me that it's called \href{https://hackage.haskell.org/package/these}{\texttt{These}} in Haskell. A version of it exists in any distributive monoidal category $(\cat{C},0,+,I,\odot)$, given by $c\vee c'\coloneqq c+(c\odot c')+c'$, and note that $\vee$-monoids are the same as $\odot$-semigroups. I learned the $\vee_S$ for $S\geq 2$, as well as other monoidal structures, from a \href{https://mathoverflow.net/questions/155939/what-other-monoidal-structures-exist-on-the-category-of-sets}{mathoverflow post}.}
 where
\begin{equation}
	A\vee_SB\coloneqq A+A\times S\times B + B.
\end{equation}
We denote the $S=1$ case simply by $\vee\coloneqq\vee_1$. When $S=0$ this is the usual coproduct, $A\vee_0B\cong A+B$, which we will treat seperately since it is very important. The other important monoidal product on $\smset$ for us is $\times$.%
\footnote{
	We sometimes denote products using juxtaposition, $AB\coloneqq A\times B$. We may also do this for polynomials $pq\coloneqq p\times q$.
}
These lead to the following symmetric monoidal products on $\poly$:
\begin{align}
\label{eqn.times}
	p\times q&\coloneqq\sum_{(I,J): p(1)\times q(1)}\yon^{p[I]+ q[J]}\\
\label{eqn.otimes}
	p\otimes q&\coloneqq\sum_{(I,J): p(1)\times q(1)}\yon^{p[I]\times q[J]}\\
\label{eqn.ovee}
	p\ovee_S q&\coloneqq\sum_{(I,J): p(1)\times q(1)}\yon^{p[I]\vee_S q[J]}
\end{align}
The first two are highly relevant: the first ($\times$) is the categorical product, and the second ($\otimes$) is the \emph{Dirichlet} product, both of which come up often in practice. Writing $\ovee\coloneqq\ovee_1$, note that there is a pullback square in $\poly$:
\begin{equation}
\begin{tikzcd}
	p\ovee q\ar[r]\ar[d]&
	p\times q\ar[d]\\
	p\otimes q\ar[r]&
	p(1)\times q(1)\ar[ul, phantom, very near end, "\lrcorner"]
\end{tikzcd}
\end{equation}
Multiplication by $\yon$ yields a normal lax monoidal functor $(\poly,1,\times)\to(\poly,\yon,\otimes)$: 
\begin{equation}\label{eqn.lax_times_tensor}
	p\yon\otimes q\yon\to (p\times q)\yon
	\qqand
	\yon\cong \yon,
\end{equation}
a lax monoidal functor $(\poly,\yon,\otimes)\to(\poly,\yon,\otimes)$:
\begin{equation}
	p\yon\otimes q\yon\to (p\otimes q)\yon
	\qqand
	\yon\to\yon^2,
\end{equation}
and a strong monoidal functor $(\poly,1,\ovee)\to(\poly,\yon,\otimes)$:
\begin{equation}
	p\yon\otimes q\yon\cong (p\ovee q)\yon
	\qqand
	\yon\cong \yon.
\end{equation}
There is a natural map%
\footnote{See also \eqref{eqn.ovee_garner}.}
\begin{equation}
(p_1\tri p_2)\ovee (q_1\tri q_2)\to(p_1\otimes q_1)\tri(p_2\ovee q_2).
\end{equation}

The Dirichlet product commutes with connected limits in either variable: for any \emph{connected} category $\cat{J}$, functor $p\colon\cat{J}\to\poly$, and polynomial $q$, the induced map
\begin{equation}\label{eqn.tensor_limits}
	\left(\lim_{j:\cat{J}}p_j\right)\otimes q
	\To{\cong}
	\lim_{j:\cat{J}}(p_j\otimes q)
\end{equation}
is an isomorphism. It commutes with all colimits in either variable: for any category $\cat{J}$, functor $p\colon\cat{J}\to\poly$, and polynomial $q$, the induced map
\begin{equation}
	\colim_{j:\cat{J}}(p_j\otimes q)
	\To{\cong}
	\left(\colim_{j:\cat{J}}p_j\right)\otimes q
\end{equation}
is an isomorphism.

The $\otimes$ and $\times$ operation together form a linearly distributive category,%
\footnote{This fact \eqref{eqn.shapiro_0}, along with \cref{eqn.shapiro_1,eqn.shapiro_2}, were discovered either by or in conjunction with Brandon Shapiro.}
\begin{equation}\label{eqn.shapiro_0}
p\times(q\otimes r)\to (p\times q)\otimes r.
\end{equation}

\chapter{Substitution product}

There is a nonsymmetric monoidal structure on $\poly$ given by composing polynomials. Its unit is $\yon$ and its monoidal product is given by the following formula:
\begin{equation}
p\tri q\coloneqq\sum_{I: p(1)}\sum_{J\colon p[I]\to q(1)}\yon^{\sum\limits_{i: p[I]}q[Ji]}
\end{equation}
If $p\to p'$ and $q\to q'$ are cartesian, then so is $p\tri q\to p'\tri q'$. If $q\to q'$ is vertical, then so is $p\tri q\to p\tri q'$. 

The functor $\ext\colon\poly\to\End(\smset)$ from \eqref{ext_functor} sends $\tri$ to composition of functors $\ext(p\tri q)\cong\ext(p)\circ\ext(q)$. Since $p\tri X$ is a set for any $p:\poly$ and $X:\smset$, we can write
\begin{equation}
\ext(p)(X)\cong p\tri X
\end{equation}
Since $\ext$ is fully faithful, $\varphi\colon p\to q$ is an isomorphism iff $(\varphi\tri X)\colon p\tri X\to q\tri X$ is a bijection for all $X:\smset$. In general, when we write $p(1)$, we really mean either $(p\tri 1):\poly$ or $\ext(p)(1):\smset$; hopefully the type is clear from context.

The monoidal structure $\tri$ is left distributive with respect to $+$ and $\times$:
\begin{align}
	0\tri q&\cong 0&&(p+p')\tri q\cong (p\tri q)+(p'\tri q)\label{eqn.comp_plus}\\
	1\tri q&\cong 1&&(p\times p')\tri q\cong (p\tri q)\times(p'\tri q)\label{eqn.comp_times}
\end{align}
In fact, $\tri$ preserves all limits in the left-variable, and it preserves connected limits in the right variable. That is, if $\cat{I}$ is any category and $\cat{J}$ is connected, then there are isomorphisms
\begin{equation}
\left(\lim_{i:\cat{I}}p_i\right)\tri p'\To{\cong}\lim_{i:\cat{I}}(p_i\tri p')
\qqand
q'\tri\left(\lim_{j:\cat{J}}q_j\right)\To{\cong}\lim_{j:\cat{J}}(q'\tri q_j)
\end{equation}
for any $p\colon\cat{I}\to\poly$, $q\colon\cat{J}\to\poly$, and $p',q':\poly$.

If $p:\poly$ is finitary (each $p[I]$ is a finite set) then for any sifted (e.g.\ filtered) category $\cat{J}$ and diagram $q\colon \cat{J}\to\poly$, the natural map
\begin{equation}\label{eqn.finitary_tri_sifted}
  \colim_{j: \cat{J}}(p\tri q_j)
  \To{\cong}
	p\tri\colim_{j: \cat{J}}q_j
\end{equation}
is an isomorphism. In fact \cref{eqn.finitary_tri_sifted} also holds for any cardinal $\kappa$ such that $p[I]<\kappa$ for all $I:p(1)$, any $\kappa$-filtered category $\cat{J}$, and any diagram $q\colon \cat{J}\to\poly$.

For any filtered category $\cat{I}$, diagram $p\colon \cat{I}\to\polycart$ of cartesian maps between polynomials $p_i$, and for any polynomial $q:\poly$, the natural map
\begin{equation}\label{eqn.cart_tri}
	\colim_{i: \cat{I}}(p_i\tri q)
	\To{\cong}
	(\colim_{i: \cat{I}}p_i)\tri q
\end{equation}
is an isomorphism.%
\footnote{The functor $(-\tri q)$ does not preserve coequalizers of cartesian maps and also does not preserve filtered colimits of arbitrary maps. A counterexample of the first is that the coequalizer of the two distinct cartesian maps $\yon^2\tto\yon^2$ is $1$, but the corresponding coequalizer of $\yon^2\tri Q\tto \yon^2\tri Q$ is isomorphic to $Q$. A counterexample of the second is that the filtered diagram $\yon^\nn\to\yon^\nn\to\cdots$ given by maps $(n\mapsto n+1)\colon\nn\to\nn$ has $1$ as its colimit, but the corresponding colimit of $\yon^\nn\tri Q$ is only $1$ for $Q=1$; for example when $Q=0$ this colimit is clearly $0$. These counterexamples were found in collaboration with Sophie Libkind.
}
 For any $p:\poly$ the operation $(p\tri\,-)$ preserves monomorphisms and epimorphisms in $\poly$. For any $q:\poly$, the operation $(-\tri q)$ preserves monomorphisms and if $q\neq 0$ then it also preserves epimorphisms.%
\footnote{The unique map $\yon\to 1 $ is an epimorphism, but $\yon\tri 0\to 1\tri 0$ is not.}

All functors $\smset\to\smset$ have a \emph{$\times$-strength}; this is expressed for any $p:\poly$ by the data of natural maps
\begin{equation}\label{eqn.strength}
	(p\tri q)\times r\to p\tri(q\times r).
\end{equation}

The monoidal structure $\tri$ is normal duoidal with $\otimes$, i.e.\ they have the same unit, $\yon$, and there is a natural transformation
\begin{equation}\label{duoidal}
	(p_1\tri p_2)\otimes(q_1\tri q_2)\too(p_1\otimes q_1)\tri(p_2\otimes q_2)
\end{equation}
satisfying the usual laws. Using $\yon$ in place of $p_1$, $p_2$, $q_1$, or $q_2$, \eqref{duoidal} induces natural maps
\begin{equation}
	[p,q]\to[r\tri p, r\tri q]
	\qqand
	[p,q]\to[p\tri r, q\tri r]
\end{equation}
and
\begin{equation}
	r\tri[p,q]\to[p,r\tri q]
	\qqand
	[p,q]\tri r\to [p,q\tri r].
\end{equation}
The identity functor $\poly\to\poly$ is lax monoidal as a functor $(\poly,\yon,\tri)\to(\poly,\yon,\otimes)$, i.e.\ for every $p,q$ a map of polynomials, which we call $\indep$:
\begin{equation}\label{eqn.indep}
	p\otimes q\To{\indep} p\tri q
\end{equation}
satisfying the usual laws. This map is derived from \eqref{duoidal} by taking $p_1\coloneqq p$, $q_1\coloneqq\yon$, $p_2\coloneqq\yon$, and $q_2\coloneqq q$. Note that if $p$ is linear or $q$ is representable, then \eqref{eqn.indep} is an isomorphism:
\begin{equation}
  A\yon\otimes q\cong A\yon\tri q
  \qqand
  p\otimes\yon^A\cong p\tri \yon^A.
\end{equation}

Since $A\yon$ and $\yon^A$ are extensionally left and right adjoint functors $\smset\to\smset$, they are left and right duals of one another in the sense that there are maps:
\begin{equation}
	\yon\to \yon^A\tri A\yon
	\qqand
	A\yon\tri\yon^A\to\yon
\end{equation}
making the triangle equations hold. If $a$ has a right dual, then $a=A\yon$ for some $A:\smset$.

The natural maps \eqref{eqn.Day}, defining $\times$ and $\otimes$ as Day convolutions of $+,\times$, can be extended from sets to polynomials using \cref{eqn.comp_plus,eqn.comp_times}. That is, we have maps:%
\footnote{Note that the map in \eqref{eqn.dayx} induces the strength \eqref{eqn.strength}, using $q_1=\yon$.}
\begin{align}
\label{eqn.dayx}
	(p_1\tri p_2)\times(q_1\tri q_2)&\to(p_1\otimes q_1)\tri(p_2\times q_2)\\
\label{eqn.day+}
	(p_1\tri p_2)\times(q_1\tri q_2)&\to(p_1\times q_1)\tri(p_2+q_2)
\intertext{It is sometimes useful to factor \eqref{eqn.day+} through a map of the form}
	(p_1\tri p_2)\times (q_1\tri q_2)&\to (p_1\times q_1) \tri\Big(p_2\times(q_1\tri q_2)+(p_1\tri p_2)\times q_2\Big).
\end{align}


\chapter{Monoidal closures}

There are closures for $\times$, $\otimes$, and $\ovee_S$ for each $S:\smset$ (\cref{eqn.times,eqn.otimes,eqn.ovee}),%
given by
\begin{align}
  q^p&\coloneqq \prod_{I: p(1)}q\tri(p[I]+\yon)\label{eqn.cart_cl}\\
  [p,q]&\coloneqq\prod_{I: p(1)}q\tri(p[I]\times\yon)\\
  \present{p,q}_S&\coloneqq\prod_{I: p(1)}q\tri(p[I]+p[I]\times S\yon+\yon)
\end{align}
These satisfy the defining universal properties:
\begin{align}
  \poly(p',q^p)&\cong\poly(p'\times p,q)\\
  \poly(p',[p,q])&\cong\poly(p'\otimes p,q)\\
  \poly(p',\present{p,q}_S)&\cong\poly(p'\ovee_S p,q)
\end{align}
The first one, $q^p$, is the cartesian closure; it is preserved by the forgetful functor $\ext\colon\poly\to\End(\smset)$. The third is not particularly interesting, as far as we know. The second one, $[p,q]$, is what we call the \emph{Dirichlet closure}; it has a very nice standard form:
\begin{equation}\label{eqn.dirichlet_hom_standard}
[p,q]\cong\sum_{\varphi:\poly(p,q)}\yon^{\sum\limits_{I: p(1)}q[\varphi_1I]}
\end{equation}
where $\varphi_1\colon p(1)\to q(1)$ is the $1$-component of the natural transformation $\varphi$. Another nice representation is
\begin{equation}
[p,q]\cong\prod_{I: p(1)}\sum_{J: q(1)}\prod_{j: q[J]}\sum_{i: p[I]}\yon.
\end{equation}

The cartesian closure satisfies all the usual arithmetic properties:
\begin{gather}\label{eqn.timeshom_properties}
	q^0\cong1,\quad
	q^{p_1+p_2}\cong (q^{p_1})\times(q^{p_2}),\quad
	1^p\cong 1,\quad
	(q_1\times q_2)^p\cong q_1^p\times q_2^p,\quad
	q^1\cong q,\quad
	q^{p_1\times p_2}\cong (q^{p_2})^{p_1}
\end{gather}
The Dirichlet closure has only some of the analogous properties:
\begin{gather}\label{eqn.dirhom_properties}
	[0,p]\cong1,\qquad
	[p_1+p_2,q]\cong [p_1,q]\times[p_2,q],\qquad
	[\yon,q]\cong q,\qquad
	[p_1\otimes p_2,q]\cong[p_1,[p_2,q]]
\end{gather}
The Dirichlet closure also satisfies the following relation with $(1,\times)$:
\begin{equation}\label{eqn.closure_products}
  [p,1]\cong 1
  \qqand
  [p,q_1\times q_2]\cong[p,q_1]\times[p,q_2].
\end{equation}
In fact, by \eqref{eqn.tensor_limits} and \eqref{eqn.closure_products}, the natural map
\begin{equation}
	[p,\lim_{i:\cat{I}}q_i]\To{\cong}\lim_{i:\cat{I}}[p,q_i]
\end{equation}
is an isomorphism for any $p:\poly$, $q\colon\cat{I}\to\poly$. By \eqref{eqn.dirichlet_hom_standard} and $-\tri 1$ is a left adjoint \eqref{eqn.adjunctions} so preserves colimits, one checks that the following is also an isomorphism
\begin{equation}
	[\colim_{i:\cat{I}}q_i,p]\To{\cong}\lim_{i:\cat{I}}[q_i,p].
\end{equation}

The following pattern holds for any monoidal closure, but we include it for $\times$ and $\otimes$ in case it's convenient:
\begin{equation}
	q_1^{p_1}\times q_2^{p_2}\to (q_1q_2)^{p_1p_2}
	\qqand
	[p_1,q_1]\otimes[p_2,q_2]\to[p_1\otimes p_2,q_1\otimes q_2]
\end{equation}

For both $\times$ and $\otimes$-closures, coproducts in the codomain reduce when the domain is representable:
\begin{equation}
	[\yon^A,p+q]\cong[\yon^A,p]+[\yon^A,q]
	\qqand
	(p+q)^{\yon^A}\cong p^{\yon^A}+q^{\yon^A}.
\end{equation}

There is a vertical epimorphism
\begin{equation}
	[p\times q,r]\to[p,r^q]
\end{equation}
which is an isomorphism iff $q$ is constant, i.e.\ $q=A$ for some $A:\smset$.

By duoidality, there are coherent natural maps
\begin{equation}
	[p_1,q_1]\tri[p_2,q_2]\to[p_1\tri p_2,q_1\tri q_2].
\end{equation}

The cartesian closure also interacts with substitution as follows:
\begin{equation}
	r\tri (q^p)\to (r\tri q)^p
\end{equation}
and this map is an isomorphism in case $p\cong\yon^A$ for some $A:\smset$. 

The functor $\yon^-\colon(\poly,1,\times)\op\to(\poly,\yon,\tri)$ is strong monoidal
\begin{equation}
	\yon\cong\yon^1
	\qqand
	\yon^p\tri\yon^q\cong\yon^{p\times q}.
\end{equation}
We will see a generalization of this, lax monoidality for any monad in place of $\yon$, in \eqref{eqn.monad_tri}.

The Dirichlet closure interacts with substitution of linear and representables as follows:
\begin{equation}\label{eqn.innerhom_adj}
  [p\tri \yon^S,q]\cong[p,q\tri S\yon]
  \qqand
  [S\yon\tri p,q]\cong[p,\yon^S\tri q]
\end{equation}
Since $[\yon,r]\cong r$ for any $r:\poly$, these induce isomorphisms:
\begin{equation}\label{eqn.closure_lin_rep}
	[\yon^S,q]\cong q\tri S\yon
	\qqand
	[S\yon,q]\cong\yon^S\tri q.
\end{equation} 

The cartesian closure interacts with coproducts and products via maps
\begin{equation}
  q^p\to(q+r)^{p+r}
  \qqand
	q^p\to(qr)^{pr},
\end{equation}
natural in $p:\poly\op$, $q:\poly$, and dinatural in $r$.

For any $a:\poly$, there is an enrichment of $\poly$ in $(\poly,\yon,\otimes)$ via
\begin{equation}
	\ul\Hom(p,q)\coloneqq[p,q^a].
\end{equation}
That is, there are natural maps $\yon\to[p,p^a]$ and $[p,q^a]\otimes[q,r^a]\to[p,r^a]$, which are appropriately unital and associative. When $a=1$ this is the usual self-enrichment.

Because of duoidality \eqref{duoidal}, $\otimes$-closure interacts with substitution:
\begin{equation}\label{eqn.shapiro_1}
	\yon\cong[\yon,\yon]
	\qqand
	[p_1,q_1]\tri[p_2,q_2]\to [p_1\tri p_2,q_1\tri q_2].
\end{equation}

Dirichlet-mapping into $\yon$ is often of interest; we have the following maps and isomorphisms:
\begin{align}
\label{eqn.homy2}
  \yon^A&\cong[A\yon,\yon]\\
\label{eqn.homy4}
  A\yon&\cong[\yon^A,\yon]\\
\label{eqn.homy1}
  [p,\yon]\times[q,\yon]&\cong[p+q,\yon]\\
\label{eqn.homy3}
  [p,\yon]+[q,\yon]&\to[pq,\yon]\\
\label{eqn.homy5}
	[p,\yon]\otimes[q,\yon]&\to[p\otimes q,\yon]\\
\label{eqn.homy6}
	[p,\yon]\tri[q,\yon]&\to[p\tri q,\yon]
\end{align}
\cref{eqn.homy1,eqn.homy3} hold for any $r:\poly$ in place of $\yon$, and \cref{eqn.homy5,eqn.homy6} hold for any $(\otimes)$-monoid in place of $\yon$; see \cref{eqn.otimes_monoids}.

\chapter{Coclosures for substitution and Dirichlet product}

The left Kan extension of a polynomial functor $p$ along another polynomial functor $q$ is again a polynomial functor, which we denote
\begin{equation}
\lens{p}{q}\coloneqq\sum_{I: p(1)}\yon^{q\tri\; (p[I])}
\end{equation}
This satisfies the following universal property of a Kan extension, i.e.\ a right-coclosure:%
\footnote{I learned the right-coclosure from Josh Meyers. I learned the in-retrospect-obvious fact that it is the same as a left Kan extension from Todd Trimble.}
\begin{equation}\label{eqn.cocl_define}
	\poly\left(\lens{p}{q},p'\right)\cong\poly\left(p,p'\tri q\right).
\end{equation}
The unit and counit of the adjunction are:
\begin{equation}\label{eqn.coclosure_unit_counit}
  p\to\lens{p}{q}\tri q
  \qqand
	\lens{p\tri q}{q}\to p
\end{equation}

The coclosure $\lens{p}{q}$ is covariant in $p$ and contravariant in $q$:
\begin{equation}
  \begin{prooftree}
  	\Hypo{p\to p'}
  	\Hypo{q\from q'}
    \Infer2{\lens{p}{q}\to\lens{p'}{q'}}
  \end{prooftree}
\end{equation}
and if $p\to p'$ is vertical or cartesian, then so is $\lens{p}{q}\to\lens{p'}{q}$, respectively. 

Note that $\lens{x}{x}$ is a comonad for any $x:\poly$,
\begin{equation}\label{eqn.lens.xx.comonad}
	\lens{x}{x}\to\yon
	\qqand
	\lens{x}{x}\to\lens{x}{x}\tri\lens{x}{x}.
\end{equation}

A formal consequence of the existence of a coclosure is that $\poly$ is enriched in $(\poly\op,\yon,\tri)$ via
\begin{equation}\label{eqn.ihom1}
	\iHom(p,q)\coloneqq\lens{p}{q}
\end{equation}
and the unit and composition are derived from \eqref{eqn.coclosure_unit_counit}. In fact, for any polynomial $x:\poly$, there is an enrichment $\iHom(p,q)\coloneqq\lens{p\tri x}{q\tri x}$. As a subcategory, $\smset\ss\poly$ is also enriched in $(\poly\op,\yon,\tri)$ via\footnote{
In both \cref{eqn.ihom1,eqn.ihom2}, the underlying set of maps is as expected:
$
\poly\op\left(\yon,\lens{p}{q}\right)=\poly\left(\lens{p}{q},\yon\right)=\poly(p,q)
$ and
$
\poly\op(\yon,A\yon^B)=\smset(A,B).
$
}
\begin{equation}\label{eqn.ihom2}
	\iHom(A,B)\coloneqq\lens{A}{B}=A\yon^B.
\end{equation}

The coclosure interacts with $\tri$ in a few different ways.
\begin{equation}\label{eqn.lens_tri}
	p\cong\lens{p}{\yon}.
	\qqand
	\lens{\lens{p}{q}}{q'}\cong\lens{p}{q'\tri q}
\end{equation}
We have a natural map
\begin{equation}\label{eqn.tricoclose}
	\lens{p\tri p'}{q\tri q'}\to \lens{p\tri\lens{p'}{q'}}{q},
\end{equation}
and in the case that either $q'=Q'\yon$ or $p=P\yon$ is linear, \eqref{eqn.tricoclose} is an isomorphism. By \eqref{eqn.lens.xx.comonad}, for any $x:\ob\poly$ there is in particular a map natural in $p,q$:
\begin{equation}
\lens{p\tri x}{q\tri x}\to\lens{p}{q}.
\end{equation}

Since $\lens{-}{q}$ is a left adjoint, it interacts with $+$ by
\begin{equation}\label{eqn.cocl_plus}
	\lens{0}{q}=0
	\qqand
	\lens{p+p'}{q}\cong\lens{p}{q}+\lens{p'}{q}.
\end{equation}
It follows that for any set $A$, we have an isomorphism
\begin{equation}\label{eqn.lens_Abase}
	\lens{Ap}{q}\cong A\times\lens{p}{q}.
\end{equation}

Using Day convolution \eqref{eqn.dayx} and the unit \eqref{eqn.coclosure_unit_counit}, coclosure interacts with $\times$ and $\otimes$ by vertical maps
\begin{equation}
 	\lens{1}{1}\cong\yon
	\qqand
  \lens{pp'}{qq'}\to\lens{p}{q}\otimes\lens{p'}{q'}.
\end{equation}
Using \eqref{eqn.day+} we also have
\begin{equation}
	\lens{1}{0}\cong 1
	\qqand
	\lens{pp'}{q+q'}\to\lens{p}{q}\lens{p'}{q'}.
\end{equation}

When $p=\yon^P$ is representable, we have
\begin{equation}\label{eqn.cocl_rep}
	\lens{\yon^P}{q}\cong\yon^{q\tri P},
\end{equation}
so we obtain \eqref{eqn.cocl_define} from \cref{eqn.cocl_plus,eqn.cocl_rep}. 

The coclosure interacts with $\otimes$ via vertical maps
\begin{equation}\label{eqn.coclosure_tensor}
	\lens{\yon}{\yon}\cong\yon
	\qqand
	\lens{p_1\otimes p_2}{q_1\otimes q_2}\to\lens{p_1}{q_1}\otimes\lens{p_2}{q_2}
\end{equation}
The latter map in \eqref{eqn.coclosure_tensor} is an isomorphism iff any of the following hold:
\begin{itemize}
	\item $q_1$ and $q_2$ are both constant,
	\item $q_1$ and $q_2$ are both linear, or
	\item $p_1$ and $p_2$ are both linear.
\end{itemize}
Here are two examples of this. First, for any sets $A,B,C:\smset$ one obtains maps
\begin{equation}
	\lens{A}{B}\otimes\lens{B}{C}\to\lens{A}{C}.
\end{equation}
using the inverse of \eqref{eqn.coclosure_tensor} and the strength \eqref{eqn.strength}. And second, when used together with \eqref{eqn.cocl_rep}, we obtain a map for any $p,q:\poly$,
\begin{equation}
	\lens{p}{q}=\lens{\yon\otimes p}{q\otimes\yon}\to\yon^{q(1)}\otimes p
\end{equation}
which is an isomorphism if either $p=P\yon$ or $q=Q\yon$ is linear.

The coclosure interacts with the $\otimes$-closure (Dirichlet-hom) by maps
\begin{equation}
\label{eqn.cocl_closure1}
	\left[\lens{p}{q_2},q_1\right]\to[p,q_1\tri q_2]
\end{equation}
\begin{equation}\label{eqn.shapiro_2}
	\lens{p_1\otimes p_2}{q}\to\lens{p_1}{\left[p_2,q\right]}
\end{equation}
The map \eqref{eqn.cocl_closure1} is vertical. It is cartesian (and hence an isomorphism) if $q_2$ is linear or $q_1$ is constant. The map \eqref{eqn.shapiro_2} is an isomorphism if $p_2$ is representable. For example, there is a map of the form
\begin{equation}
  p_1\otimes p_2\to\lens{p_1}{[p_2,\yon]}.
\end{equation}

For any sets $S,B:\smset\op$ and $A,T:\smset$, we have an isomorphism
\begin{equation}
	\left[\lens{S}{T},\lens{[B,A]}{B}\right]
	\cong
	\left[\lens{SB}{1},\lens{AT}{1}\right]
\end{equation}
For example, it's useful in the theory of Moore and Mealy machines to know that $[S\yon^T,\yon]\cong[S\yon,T\yon]$.

There is also a natural isomorphism $[p,\yon]\cong\Gamma(p)\times\lens{\yon}{p}$, where $\Gamma(p)=\poly(p,\yon)$; in particular there is a natural map
\begin{equation}
	[p,\yon]\to\lens{\yon}{p}.
\end{equation}
In fact, for any $I:\smset$ and $p\colon I\to\poly$, there is a natural vertical map
\begin{equation}
	\left[\bigotimes_{i:I}p_i,\,\yon\right]\too\bigotimes_{i:I}\left[p_i,\,\bigotimes_{j\neq i}\lens{\yon}{p_j}\right],
\end{equation}
e.g.\ $[p\otimes q,\yon]\to[p,\yon^{q(1)}]\otimes[q,\yon^{p(1)}]$ and an isomorphism
\begin{equation}
	[p\otimes q,\yon]\cong\left[\lens{p}{\yon^{q(1)}},\yon\right]\otimes\left[\lens{q}{\yon^{p(1)}},\yon\right].
\end{equation}

The coclosure interacts with cartesian closure by way of a natural map
\begin{equation}
	\lens{p\times q}{r}\to\lens{p}{r^q}.
\end{equation}
This map is almost never an isomorphism (it is if and only if $p=0$ or $q=1$). However, for any $A:\smset$, it is always a ``$\yon^A$-local equivalence'', i.e.\ it induces a bijection
\begin{equation}
	\poly\left(\lens{p}{r^q},\yon^A\right)
	\cong
	\poly\left(\lens{pq}{r},\yon^A\right).
\end{equation}

All functors between cartesian categories are oplax monoidal; in the case of the functor $\lens{p}{-}$ for any $p:\poly$, the counit and comultiplication maps factor through isomorphisms:
\begin{equation}\label{eqn.lens_plus_upstairs_pullback}
	\lens{p}{0}\To{\cong}p(1)
	\qqand
	\lens{p}{q_1+q_2}\To{\cong}\lens{p}{q_1}\times_{p(1)}\lens{p}{q_2}
\end{equation}

For any $p$, we have an isomorphism (left), and it follows from \eqref{eqn.lens_plus_upstairs_pullback} that for any $B:\smset$ there are also isomorphisms (middle and right)
\begin{equation}
  \lens{p}{1}\cong p(1)\yon
  \qqand
  \lens{p}{B}\cong p(1)\yon^B
  \qqand
  \lens{p}{B+\yon}\cong\yon^B\times p
\end{equation}

For any symmetric monoidal product $(I,\cdot)$ on $\smset$, let $(\yon^I,\odot)$ be the associated symmetric monoidal product on $\poly$ (see \cref{eqn.day_all}). Then for any $p,q:\poly$ there is an isomorphism
\begin{equation}
	p\odot q\cong\sum_{I:p(1)}\lens{q}{\yon\odot p[I]}.
\end{equation}
From this and \eqref{eqn.day_distributes} one can obtain the monoidal closure formula, e.g.\ \eqref{eqn.cart_cl}.

For any $A:\smset$ and $p:\poly$ we have isomorphisms
\begin{equation}\label{eqn.left_adjoints}
  \lens{p}{A\yon}\cong 
  p\tri\yon^A\cong 
  p\otimes\yon^A
  \qqand
  A\yon\otimes\lens{p}{q}\cong
  A\yon\tri\lens{p}{q}\cong
	\lens{A\yon\tri p}{q}
\end{equation}
\cref{eqn.left_adjoints} generalize to bicomodules---namely $A\yon$ and $\yon^A$ can be replaced by any pair of left and right adjoint bicomodules---even though that is beyond the scope of this document. Indeed, quite a few of the structures in this document generalize to the bicomodule setting.

There is also an \emph{indexed} left $\tri$-coclosure. That is, for any function $f\colon p(1)\to q(1)$, define
\begin{equation}
	p\cocl{f}q\coloneqq \sum_{I: p(1)}q[fI]\,\yon^{p[I]}.
\end{equation}
This satisfies the following indexed-adjunction formula:%
\footnote{Note that the indexed adjunction \eqref{indexed_adjunction} is not natural in $q:\poly$, but it is natural in $q:\polycart$.}
\begin{equation}\label{indexed_adjunction}
	\poly(p,q\tri r)\cong\sum_{f\colon p(1)\to q(1)}\poly(p\cocl{f}q,r)
\end{equation}
Given $\varphi\colon p\to q\tri r$, we denote its image under this isomorphism by $(\varphi.1,\varphi.2)$, where $\varphi.1\colon p(1)\to q(1)$ and $\varphi.2\colon (p\cocl{\varphi.1}q)\to r$.

The indexed coclosure is a very well-behaved structure.
\begin{align}
	(p\cocl{!}\yon^A)&\cong Ap\\
	(p+p')\cocl{(f,f')}q&\cong(p\cocl{f}q)+(p'\cocl{f'}q)\\
	(p+p')\cocl{f+f'}(q+q')&\cong(p\cocl{f}q)+(p'\cocl{f'}q')\\
	p\cocl{(f,f')}(q\times q')&\cong(p\cocl{f}q)+(p\cocl{f'}q')\\
	(p\times p')\cocl{\pi_1\then f}q&\cong(p\cocl{f}q)\times p'\\
	(p\cocl{f}q)\times p'&\cong(p\times p')\cocl{(p\times!)\then f}q\\
	(p\cocl{f}q)\tri p'&\cong(p\tri p')\cocl{(p\tri!)\then f}q\\
	p\cocl{f}(q\tri r)&\cong(p\cocl{f.1}q)\cocl{f.2}r\\
	(p\otimes p')\cocl{f\otimes f'}(q\otimes q')&\cong(p\cocl{f}q)\otimes(p'\cocl{f'}q')\\
	\lens{p}{r}\cocl{f}q&\cong\lens{p\cocl{f}q}{r}
\end{align}

Given $p,q:\poly$ and map $f\colon p(1)\to q(1)$, there is a natural map
\begin{equation}
	\lens{p}{p\cocl{f}q}\to q.
\end{equation}
For example when $p=\yon$, this is the cartesian map picking out the representable summand of $q$ at position $f$.

Using $\frown$ we can define for any $p:\poly$ an object $p_*:\poly$ as follows.
\begin{equation}\label{eqn.p_star}
p_*\coloneqq p\cocl{\id}p=\sum_{I: p(1)}p[I]\yon^{p[I]}
\end{equation}
The assignment $p\mapsto p_*$ is not functorial in $\poly$, though it is functorial in the cartesian morphisms
\begin{equation}
	(p\mapsto p_*)\colon\polycart\to\polycart
\end{equation}
The bundle representation of a polynomial $p$ is $p_*(1)\to p(1)$.

In fact, the $p\mapsto p_*$ construction is strong monoidal with respect to $+$, $\otimes$, and hence $\OR$:%
\footnote{Although there is a map of polynomials $\dot{p}\otimes\dot{q}\to\dot{(p\otimes q)}$, it is not cartesian, and hence does not constitute any sort of lax monoidal structure on $p\mapsto\dot{p}$.}
\begin{align}
	0_*\cong 0
	&\qqand
	(p+q)_*\cong p_*+q_*\\
	\yon_*\cong\yon
	&\qqand
	(p\otimes q)_*\cong p_*\otimes q_*\\
	0_*\cong0
	&\qqand
	(p\OR q)_*\cong p_*\OR q_*
\end{align}
This construction extends to a comonad $p_*\to p$ and $p_*\to (p_*)_*$ on $\polycart$. 

For each $p:\poly$, its pointing $p_*$ is a $\tri$-comonoid
\begin{equation}\label{eqn.pointing_comonad}
	\epsilon\colon p_*\to\yon
	\qqand
	\delta\colon p_*\to (p_*)_*.
\end{equation}
The map $\epsilon$ sends $(i,d_0)\mapsto d_0$ and the map $\delta$ sends $(i,d_0)\mapsto((i,d_0), d\mapsto ((i,d), d'\mapsto d'))$. So as a category $p_*$ has objects $\sum_{i:p(1)}p[i]$, and for each object $(i,d_0)$, a map out is an element $d:p[i]$; its codomain is $(i,d)$, and composing is just taking the most recent pointing.

The composition map $[p,q]\otimes [q,r]\to [p,r]$ factors via $\epsilon\colon q_*\to \yon$ through a map
\begin{equation}\label{eqn.eval.factor}
	[p,q]\otimes [q,r]\to [p,r]\tri q_*\;,
\end{equation}
For fixed $q$, \cref{eqn.eval.factor} is functorial in $p:\poly\op$ and $r:\poly$, and its mate $[p,q]\to[[q,r],[p,r]\tri q_*]$ is functorial in $q:\polycart$ for fixed $p,r$. In particular, there is a map
\begin{equation}
p\otimes[p,\yon]\to p_*.
\end{equation}

The operation $(p\mapsto p_*)$ is related to the derivative $p\mapsto\dot{p}$, where
\begin{equation}
\dot{p}\coloneqq\sum_{I:p(1)}\sum_{i:p[I]}\yon^{p[I]-\{i\}},
\end{equation}
indeed $p_*\cong\dot{p}\yon$. The derivative operation is an adjoint in $\polycart$:
\begin{equation}
  \polycart(q\yon,p)\cong
  \polycart(q,\dot{p}).
\end{equation}
There are maps $p\to p\tri\dot{p}$, $p\to p\tri\dot{p}\tri\ddot{p}$, etc., all natural in $\polycart$. 

The product and chain rules from calculus hold
\begin{equation}\label{eqn.calculus}
	\dot{(pq)}\cong \dot{p}q+p\dot{q}
	\qqand
	\dot{(p\tri q)}\cong(\dot{p}\tri q)\dot{q}.
\end{equation}
and they extend to isomorphisms 
\begin{equation}\label{eqn.calculus}
	(pq)_*\cong p_*q+pq_*
	\qqand
	(p\tri q)_*\cong(\dot{p}\tri q)q_*.
\end{equation}

Given a map $\varphi\colon p\to p'$ and a function $g\colon q(1)\to q'(1)$, there is a cartesian map
\begin{equation}
	(p\cocl{\varphi\tri1}p')\tri(q\cocl{g}q')\to
	(p\tri q)\cocl{\varphi\tri g}(p'\tri q'),
\end{equation}
natural in $p,q:\poly$ and $p',q':\polycart$. In particular, there is a cartesian map
\begin{equation}
	p_*\tri q_*\to(p\tri q)_*
\end{equation}
which is part of a normal lax monoidal functor $-_*\colon(\poly,\yon,\tri)\to(\poly,\yon,\tri)$.

For any polynomial $p$, there is a map
\begin{align}
	(p_*)_*&\to\yon+1=\sum_{i:\{0,1\}}\yon^i\\
	(i,d_1,d_2)&\mapsto
	\begin{cases}
		(1,d_1)&\text{ if } d_1=d_2 \\
		(0,!)&\text{ if }d_1\neq d_2
	\end{cases}
\end{align}
The vertical-cartesian factorization of $(p_*)_*\to\yon+1$ was called the \emph{identity type} for $p$ in \cite{aberle2024polynomial}; it is
\begin{equation}
	p^=\coloneqq\sum_{i:p(1)}\sum_{d_1,d_2:p[i]}\yon^{\begin{cases}1&\text{ if } d_1=d_2 \\0&\text{ if }d_1\neq d_2\end{cases}}
\end{equation}
There is a cartesian map
\begin{equation}
	p\to (p\otimes p)\tri p^=
\end{equation}
Note that $p^=$ is affine 

The assignment $p\mapsto p_*$ is also functorial as a map $\poly\to\Span(\poly)$:
\begin{align}
	\cocl{}\colon\poly&\to\Span(\poly)\\
	p&\mapsto p_*\\
	(p\To{\varphi}q)&\mapsto\big( p_*\from(p\cocl{\varphi_1}q)\to q_*\big)
\end{align}\goodbreak
That is, for any $p\To{\varphi}q\To{\psi}r$, there is an isomorphism
\begin{equation}
	(p\cocl{\varphi_1}q)\times_{ q_*}(q\cocl{\psi_1}r)\cong p\cocl{(\varphi\then\psi)_1}r
\end{equation}
This functor is strong monoidal with respect to both $+$ and $\otimes$. One may think of it as representing the bundle view of $\poly$. Indeed, for any $p:\poly$ we have a counit map $\epsilon_p\colon p_*\to p$, and given $\varphi\colon p\to q$, there is an induced span
\begin{equation}
\begin{tikzcd}
	 p_*\ar[d]&p\cocl{\varphi_1}q\ar[d]\ar[l]\ar[r]& q_*\ar[d]\\
	p\ar[r, equal]&p\ar[r,"\varphi"']&q
\end{tikzcd}
\end{equation}
and evaluating at $1$ returns the usual bundle picture, since $(p\cocl{\varphi_1}q)(1)\cong p(1)\times_{q(1)}q_*(1)$. Since cartesian maps in $\poly$ are exponentiable%
\footnote{I believe the cartesian maps are precisely the exponentiable maps in $\poly$.}
the maps $p_*\to p$ can be treated as $\poly$-polynomials in the sense of \cite{weber2015polynomials}, thus giving an embedding of $\poly$ into what might be denoted $\poly_\poly$.

There is a normal lax duoidal functor, which we call the \emph{Kan residual}
\begin{equation}
	\at{-}{-}\colon(\polycart,1,\times,\yon,\otimes,\tri)\times(\poly,1,\times,\yon,\otimes,\tri)\to(\polycart,1,\times,\yon,\otimes,\tri)
\end{equation}
defined by
\begin{equation}
		\at{p}{q}\coloneqq p\frown\lens{p}{q}.
\end{equation}
The formula for it can be written many ways%
\footnote{The notation $\at{p}{q}$ is supposed to very lightly evoke the notation from logic where $\phi[x/y]$ means substitute $x$'s in place of the $y$'s found in $\phi$. The function $d\colon q[J]\to p[I]$ in \eqref{eqn.at_main} is what let's us substitute $p$-directions in place of $q$-directions. See also \eqref{eqn.at_add2}, \eqref{eqn.at_p_1}, and \eqref{eqn.at_q_coalg_poly_over_q} for other justifications of this notation.}
\begin{align}
	\at{p}{q}&\cong\label{eqn.at_main}
	\sum_{I:p(1)}\sum_{J:q(1)}\sum_{d\colon q[J]\to p[I]}\yon^{p[I]}\\&\cong
	\sum_{I:p(1)}(q\tri p[I])\yon^{p[I]}\\&\cong
	\sum_{J:q(1)}p_{\underbrace{*\cdots*}_{q[J]\text{ times}}}
\end{align}
Its defining property, which drops any mention of cartesianness, is
\begin{equation}
	\poly(\at{p}{q}, r)\cong\poly_{/p(1)}\left(\lens{p}{r},\lens{p}{q}\right).
\end{equation}
The normal lax duoidality says that $\yon\cong\at{\yon}{\yon}$ and that we have coherent cartesian maps
\begin{align}
\label{eqn.atplus}
	\at{p}{q}+\at{p'}{q'}&\tocart\at{(p+p')}{(q+ q')}\\
\label{eqn.attimes}
	\at{p}{q}\times\at{p'}{q'}&\tocart\at{(p\times p')}{(q\times q')}\\
\label{eqn.atotimes}
	\at{p}{q}\otimes\at{p'}{q'}&\tocart\at{(p\otimes p')}{( q\otimes q')}\\
\label{eqn.attri}
	\at{p}{q}\tri\at{p'}{q'}&\tocart\at{(p\tri p')}{( q\tri q')},
\intertext{and there is also a similar-looking cartesian map%
\footnote{The $\upup$ map will be defined in \eqref{eqn.upup_formula}.}}
	\at{p}{q}\upup\at{p'}{q'}&\tocart\at{(p\upup p')}{(q\tri q')}.
\end{align}

For any $p,q$, the Kan residual $\at{p}{q}$ forms the apex of a span of natural maps, one of which is cartesian:
\begin{equation}\label{eqn.KR_span}
  \at{p}{q}\tocart p
  \qqand
  \at{p}{q}\to q.
\end{equation}
For any $p,q$, the functors $\at{p}{-},\at{-}{q}\colon\poly\to\polycart$ are cocartesian monoidal:
\begin{align}
\label{eqn.at_add}
\at{p}{0}&\cong 0
&
\at{p}{(q_1+q_2)}&\cong \at{p}{q_1}+\at{p}{q_2},\\
\at{0}{q}&\cong 0
&\label{eqn.at_add2}
\at{(p_1+p_2)}{q}&\cong \at{p_1}{q}+\at{p_2}{q}.
\end{align}
And the induced functor 
\begin{equation}\label{eqn.at_strong_monoidal_times}
(\poly,1,\times)\To{q\mapsto\at{-}{q}}(\End(\polycart),\id,\circ)
\end{equation}
is strong monoidal. That is, for any $p,q,q'$ we have isomorphisms
\begin{equation}\label{eqn.at_p_1}
	\at{p}{1}\cong p
	\qqand
	\at{p}{( q\times q')}\cong \at{\at{p}{q}}{q'}.
\end{equation}
Moreover, the functor in \eqref{eqn.at_strong_monoidal_times} in fact sends each $q$ to a colax duoidal endofunctor:
\begin{gather}
	\at{-}{q}\colon(\polycart,\yon,\otimes, \tri)\To{\tn{colax}}(\polycart,\yon,\otimes, \tri)\\
	\at{(p\tri p')}{q}\tocart\at{p}{q}\tri\at{p'}{q}
	\qqand
	\at{(p\otimes p')}{q}\tocart\at{p}{q}\otimes\at{p'}{q}\\
	\at{\yon}{q}\tocart\yon
\end{gather}
and the two maps $\at{\big((p_1\tri p_2)\otimes(p'_1\tri p'_2)\big)}{q}\to\big(\at{p_1}{q}\otimes\at{p'_1}{q}\big)\tri\big(\at{p_2}{q}\otimes\at{p'_2}{q}\big)$ agree.%
\footnote{Note, however, that the only $(\otimes)$- or $(\tri)$-comonoids in $\polycart$ are linear, i.e.\ of the form $A\yon$ for some $A:\smset$.}
For any $q:\poly$, the functor $\at{-}{q}\colon\polycart\to\polycart$ is a comonad; its category of coalgebras is
\begin{equation}\label{eqn.at_q_coalg_poly_over_q}
	\at{-}{q}\coalg\cong\poly_{/q}.
\end{equation}
Given a coalgebra $p\to\at{p}{q}$, one composes with $\at{p}{q}\to q$ from \eqref{eqn.KR_span} to obtain the map $p\to q$.

By \eqref{eqn.at_add} and \eqref{eqn.at_p_1}, for any $A:\smset$, we have
\begin{equation}
\at{p}{A}\cong Ap.
\end{equation}
The obvious map $\at{p}{(q\times q')}\to\at{p}{q}\times\at{p}{q}'$ factors through an isomorphism
\begin{equation}
	\at{p}{(q\times q')}\cong\at{p}{q}\times_{p(1)}\at{p}{q}'.
\end{equation}
There is also a natural cartesian map
\begin{equation}
\at{p}{q}\tocart\at{p}{(p\times q)}.
\end{equation}

There are cartesian natural maps%
\footnote{The $\upup$ map will be defined in \eqref{eqn.upup_formula}.}
\begin{align}
  p\tri\at{q}{r}&\tocart \at{(p\tri q)}{(p\tri r)}\\
  p\otimes\at{q}{r}&\tocart \at{(p\otimes q)}{(p\otimes r)}\\
  p\upup\at{q}{r}&\tocart\at{(p\upup q)}{(p\upup r)}
\end{align}

The category of polynomials is enriched in $\sspan(\poly)$:
\begin{equation}
  p\to\at{p}{p}
  \qqand
	\at{p}{q}\times_q\at{q}{r}\cong\at{\at{p}{q}}{r}\too \at{p}{r}
\end{equation}
and this enriched category is monoidal under $(+,\times,\otimes,\tri)$, by \eqref{eqn.atplus}, \eqref{eqn.attimes}, \eqref{eqn.atotimes}, \eqref{eqn.attri}. In particular $p\from\at{p}{p}\to p$ is an internal category object in $\poly$, natural in $p:\polycart$, and each $\at{p}{q}$ is an internal profunctor. 

There is a natural map
\begin{equation}
	\at{p}{q}\to(p\otimes q)\tri(\yon+1)
\end{equation}
sending $(I,J,\ol{j})\mapsto(I,J,(i,j)\mapsto \yon^{i=\ol{i}(j)}$,%
\footnote{
Here, $I:p(1)$, $J:q(1)$, $\ol{i}\colon q[J]\to p[I]$, $i:p[I]$, $j:q[J]$.
}
and it plays well with $+$, $\times$, $\otimes$, and $\tri$.

There is an indexed coclosure for $\otimes$.%
\footnote{I learned about this indexed coclosure $\hyper{}$ for $\otimes$ from Nelson Niu.}
For any function $f\colon p(1)\to q(1)$, define
\begin{equation}
p \hyper{f} q\coloneqq\sum_{I: p(1)}\yon^{\left(p[I]^{q[fI]}\right)}
\end{equation}
This satisfies the following indexed-adjunction formula:
\begin{equation}
	\poly(p,q\otimes r)\cong\sum_{f\colon p(1)\to q(1)}\poly(p\hyper{f}q,r)
\end{equation}
It also satisfies the following:
\begin{align}
	(p_1+p_2)\hyper{(f_1,f_2)}q&\cong(p_1\hyper{f_1}q)+(p_2\hyper{f_2}q)\\
	p\hyper{(f_1,f_2)}(q_1\otimes q_2)&\cong (p\hyper{f_1}q_1)\hyper{f_2}q_2
\end{align}
and for any $f\colon p(1)\to (q_1\tri q_2)(1)$, there is a natural map coming from duoidality \eqref{duoidal}:
\begin{equation}
	(p\cocl{f_1}q_1)\hyper{f_2}q_2\to
	p\hyper{f}(q_1\tri q_2).
\end{equation}
For any $f\colon p(1)\to p'(1)$ there is a natural map
\begin{equation}
\lens{p}{q\tri p'}\to\lens{p\hyper{f}p'}{q}.
\end{equation}

For any polynomial $p$, let $\lin{p}\coloneqq p(1)\yon$. For any $p,q:\poly$ there is an isomorphism
\begin{equation}
	p\cocl{}\left(\lens{\lin{p}}{q}\hyper{}p\right)\cong p\tri\lin{q}
\end{equation}
where the unwritten indices of $\cocl{}$ and $\hyper{}$ are both the identity on $p(1)$.

For any $q:\poly$ and $A:\smset$, let $q\coalg[A]:\smset$ denote the set $\smset(A,q(A))$ of $q$-coalgebra structures on $A$. For a polynomial $p$, let $q\coalg_p\coloneqq\sum_{I:p(1)}\yon^{q\coalg[p[I]]}$. Then there is an isomorphism%
\footnote{As usual, the unwritten index of $\hyper{}$ in \cref{eqn.hyper_coalg} is the identity on $p(1)$.}
\begin{equation}\label{eqn.hyper_coalg}
	q\coalg_p
	\cong
	\lens{p}{q}\hyper{}p
\end{equation}
For example, $\yon^{q\coalg[A]}\cong\lens{\yon^A}{q}\hyper{}\yon^A$.

\chapter{Other monoidal structures and actions}

We can follow the pattern we learned from Garner to take known monoidal products $\star$ that happen to distribute over $+$, and obtain new monoidal products
\begin{equation}
p\garner{\star} q\coloneqq p+(p\star q) + q.
\end{equation}
with unit 0.  monoidal products on $\poly$ from any known ones that distribute over $+$. For example, there is a symmetric monoidal structure $\OR\coloneqq\garner\otimes$ on $\poly$ with unit $0$ and product given by
\begin{equation}\label{eqn.vee}
  p\OR q\coloneqq p+(p\otimes q)+q.
\end{equation}
The functor $(p\mapsto p+\yon)$ is strong monoidal $(\poly,0,\OR)\to(\poly,\yon,\otimes)$, i.e.\ there is a natural isomorphism
\begin{equation}
	(p+\yon)\otimes(q+\yon)\cong (p\OR q)+\yon.
\end{equation}
Note that although there is an inclusion $p\otimes q\to p\OR q$, there is no map $\yon\to0$, so the identity on $\poly$ is neither lax nor oplax as a map $(\poly,0,\OR)\to(\poly,\yon,\otimes)$.

The identity functor $\poly\to\poly$ has a lax monoidal structure,
\begin{equation}
	p+q\to p\OR q.
\end{equation}
There are two duoidal structures for $\OR$, one with $\tri$ and one with $\otimes$:
\begin{align*}
	(p_1\tri p_2)\OR(q_1\tri q_2)&\to(p_1\OR q_1)\tri(p_2\OR q_2)\\
	(p_1\otimes p_2)\OR(q_1\otimes q_2)&\to(p_1\OR q_1)\otimes(p_2\OR q_2)
\end{align*}

Another very similar symmetric monoidal structure we will consider has unit $0$ and product given by
\begin{equation}\label{eqn.garner_ovee}
	p\garner{\ovee} q\coloneqq p+p\ovee q+q=p+\left(\sum_{I:p(1)}\sum_{J:q(1)}\yon^{p[I]+p[I]q[J]+q[J]}\right)+q
\end{equation}
where $\ovee$ was defined in \eqref{eqn.ovee}. There is a duoidal structure
\begin{equation}\label{eqn.ovee_garner}
	(p_1\tri p_2)\garner{\ovee} (q_1\tri q_2)\to(p_1\garner{\ovee} q_1)\tri(p_2\garner{\ovee} q_2).
\end{equation}

Of course, a third symmetric monoidal structure on $\poly$ is $\garner{\times}$, though we know no use for it.

We know of two more monoidal products $(\yon,\dagger)$ and $(\yon, \ddagger)$ from Nelson Niu, who cites de Paiva's notion of \emph{cross product} as inspiration for $\dagger$:
\begin{align}
	p\dagger q\coloneqq\sum_{I: p(1)}\sum_{J: q(1)}\prod_{i: p[I]}\prod_{j\colon p(1)\to q[J]}\yon\\
	p\ddagger q\coloneqq\sum_{I: p(1)}\sum_{J: q(1)}\prod_{i\colon q(1)\to p[I]}\prod_{j\colon p(1)\to q[J]}\yon
\end{align}
E.g.\ in the case of monomials we have $I\yon^A\dagger J\yon^B\cong IJ\yon^{AB^I}$ and $I\yon^A\ddagger J\yon^B\cong IJ\yon^{A^JB^I}$. There are natural maps
\begin{equation}
	p\otimes q\to p\dagger q
	\qqand
	p\dagger q\to p\ddagger q
\end{equation}
given by $(I,J)\mapsto ((I,J),(i,j)\mapsto(i,j(I))$ and $(I,J)\mapsto((I,J),(i,j)\mapsto(i(J),j)$, respectively.

\bigskip
Monoidal products on $\poly$ give rise to actions of $\poly$ on itself; we now turn to an action of $\polycart$ on $\poly$.
In \cite{aberle2024polynomial}, the authors (myself and Aberl\'e) define an action of $(\polycart,\yon,\tri)$
\begin{equation}
  \upup\colon\polycart\times\poly\to\poly
\end{equation}
given by
\begin{equation}\label{eqn.upup_formula}
  p\upup q\coloneqq\sum_{I:p(1)}\sum_{J\colon p[I]\to q(1)}\yon^{\prod_{d:p[I]}q[J(d)]}\;.
\end{equation}
Indeed, there are unitor and actor isomorphisms of the form
\begin{equation}
	\yon\upup q\cong q
	\qqand
	(p_1\tri p_2)\upup q\cong p_1\upup(p_2\upup q).
\end{equation}
It is also strong monoidal with respect to $(0,+)$:
\begin{equation}
	0\upup q\cong 0
	\qqand
	(p_1+p_2)\upup q\cong(p_1\upup q) + (p_2\upup q).
\end{equation}
There is an additional colax structure:
\begin{equation}
	p\upup\yon\to\yon
	\qqand
	p\upup(q_1\tri q_2)\too(p\upup q_1)\tri(p\upup q_2).
\end{equation}
All told, the above can be summarized as a strong rig monoidal functor
\begin{equation}\label{eqn.upup_summary}
\big(\polycart,0,+,\yon,\tri\big)\too\big(\Cat{ColaxEnd}(\poly,\yon,\tri),0,+,\id,\circ\big).
\end{equation}
There is also a natural map
\begin{equation}
	p\tri q\to(p\upup q)\otimes p
\end{equation}
given by $(I.J)\mapsto(I.J,I),(e,d)\mapsto(d,e(d))$.%
\footnote{Here $I:p(1)$, $J:p[I]\to q(1)$, $e:\prod_{d:p[I]}q[J(d)]$, $d:p[I]$.}

Often the $\upup$ construction is used in defining maps $p\upup q\to q$;
\footnote{One can think of a map $p\upup q\to q$ as a structure where the positions of $q$ have $p[-]$-ary joins and the directions of $q$ constitute a presheaf on this ``lattice''. For example, if $A$ is a suplattice, $\cat{F}\colon A\op\to\smset$ is a presheaf on it, and $q\coloneqq\sum_{a:A}\yon^{\cat{F}(a)}$, then there is a map $p\upup q\to q$ given by
\[
  (I,Q)\mapsto
  \left(
  	\bigvee_{d:p[I]}Q(d), \gamma\mapsto d\mapsto \gamma\big|_d
  \right)
\]
where $I:p(1)$, $Q\colon p[I]\to q(1)$, $\gamma:\cat{F}(\bigvee_{d:p[I]}Q(d))$, and $\gamma\big|_{d}$ is the restriction of the ``global'' section $\gamma$ to $d$.
}
 for example, this comes up in defining distributive laws of various sorts; we won't go into it here, but we note that for any $r:\poly$, there is a bijection
\begin{equation}
\poly(p\upup q,r)\cong\{\varphi\colon p\tri q\to r\otimes p\mid\varphi_1(I,J)=(K,I)\text{ and }\varphi^\sharp_{(I,J)}(f,d)=(d,e)\}.
\end{equation}

We note in passing that for any $p:\poly$ there is an isomorphism
\begin{equation}
p\upup\yon\cong\lens{p}{1}
\end{equation}
giving two different constructions of the ``linearization'' comonad $p\mapsto (p\tri1)\yon$; see also \eqref{eqn.adjunctions}.

There is a monoidal functor 
\begin{equation}
	\cat{G}_-\colon\Big((\polycart)\op,0,+,1,\times,\yon,\otimes\Big)\to(\grp,\{e\},\times)
\end{equation}
with $\cat{G}_p\coloneqq\prod_{I:p(1)}\Aut(p[I])$. It is strong monoidal in $0,+)$, i.e.
\begin{equation}
	\{e\}\cong\cat{G}_0
	\qqand
	\cat{G}_p\times\cat{G}_q\cong\cat{G}_{p+q},
\end{equation}
and, since it is very simple on affines like $1$ and $\yon$,
\begin{equation}
\{e\}=\cat{G}_{A\yon+B},
\end{equation}
it is normal lax monoidal in $\times,\otimes$:
\begin{align}
	\cat{G}_{p}\times\cat{G}_q\to\cat{G}_{pq}
	\qqand
	\cat{G}_{p}\times\cat{G}_q\to\cat{G}_{p\otimes q}.	
\end{align}
The group $\cat{G}_p$ acts freely and transitively on $\polycart(p,u)$, where $u\coloneqq\sum_{A:\mathbb{U}}\yon^{[A]}$ is the polynomial associated to any big enough universe $\mathbb{U}$. That is, $\polycart(p,u)$ is a $\cat{G}_p$-torsor.

There are actions of $\cat{G}_p$ and $\cat{G}_q$ on $p\tri q$. Considering $p$ and $q$ as bundles, there is a bijection
\begin{equation}
	\Cat{Dir}(p,q)\coloneqq\left\{
	\begin{tikzcd}[ampersand replacement=\&]
		p_*(1)\ar[r]\ar[d]\&q_*(1)\ar[d]\\
		p(1)\ar[r]\&q(1)
	\end{tikzcd}
	\right\}
	\cong
	\polycart(p,q\tri u)/\cat{G}_u.
\end{equation}

\chapter{Polynomial monoids and comonoids}\label{chap.monoids.comonoids}

We can consider monoids and comonoids in other monoidal structures on $\smset$. 

\paragraph{$(+)$-monoids and comonoids.}
Every polynomial is a $(+)$-monoid in a unique way, and the only object that is a $(+)$-comonoid is $0$.

\paragraph{$(\times)$-monoids and comonoids.}
Every polynomial is a $(\times)$-comonoid in a unique way.

A $(\times)$-monoid $(m,\eta,\mu)$ can be identified with a polynomial functor $m\colon\smset\to\smset$ that factors as
\begin{equation}
	\smset\To{(m,\eta,\mu)}\Cat{Monoid}\To{U}\smset
\end{equation}
where $U$ takes a monoid to its underlying set.

Any monoid $(M,e,*)$ in $\smset$ is also a $(\times)$-monoid in $\poly$. Every commutative $(\times)$-monoid $m$ in $\poly$ arises in this way, i.e.\ $m=m(1)$ is constant.

The free $(\times)$-monoid is $\uu$, defined as follows
\begin{equation}\label{eqn.list}
\uu\coloneqq\sum_{N:\nn}\yon^N\cong 1+\yon+\yon^2+\yon^3+\cdots
\end{equation}
The unit $1\to\uu$ is the empty list and $\uu\times\uu\to\uu$ is list concatenation. For any $p:\poly$ there is an isomorphism
\begin{equation}\label{eqn.unfolding_list}
\uu\tri p\cong 1+p\times(\uu\tri p)
\end{equation}
and this is the free $(\times)$-monoid on $p$.

A representable $\yon^S$ carries a monoid structure iff $S=0$, in which case it is the terminal monoid. We know from \eqref{eqn.t^c_times_monoid} that $t^c$ carries a $(\times)$-monoid structure for any polynomial monad $t$ and comonad $c$.

Suppose that $m,n$ are $(\times)$-monoids, $a$ is a $(\otimes)$-monoid, $t$ is a $(\tri)$-monoid, and $p, q$ are any polynomials; then:
\begin{align}
	mn&\quad\tn{carries the structure of a $(\times)$-monoid}\\
	m+p&\quad\tn{carries the structure of a $(\times)$-monoid}\\
	m\tri p&\quad\tn{carries the structure of a $(\times)$-monoid}\\
	[p,m]&\quad\tn{carries the structure of a $(\times)$-monoid}\\
	m^p&\quad\tn{carries the structure of a $(\times)$-monoid}\\
	a\tri m&\quad\tn{carries the structure of a $(\times)$-monoid}\\
	(mp)^p&\quad\tn{carries the structure of a $(\times)$-monoid}\\
	p^{qp}&\quad\tn{carries the structure of a $(\times)$-monoid}
\end{align}
\
\paragraph{$(\otimes)$-monoids and comonoids.}
The $(\otimes)$-monoids
\begin{equation}\label{eqn.otimes_monoids}
	\yon\To{\eta} p
	\qqand
	p\otimes p\To{\mu} p
\end{equation}
are precisely the lax monoidal polynomial functors $\smset\to\smset$, sometimes called \emph{applicatives}; indeed, by Day convolution we have $(p\tri A)\times (p\tri B)\To{\text{Day}} (p\otimes p)\tri(AB)\To{\mu}p(A\times B)$. Any $(\tri)$-monoid has an underlying $(\otimes)$-monoid structure.

The free $(\otimes)$-monoid on $p$ is 
\begin{equation}
	\freecol_p\coloneqq\sum_{N:\nn}p^{\otimes N}=\sum_{N:\nn}\sum_{i\colon N\to p(1)}\yon^{\prod_{n:N}p[i(n)]}
\end{equation}
The functor $\freecol=\freecol_-\colon\poly\to\poly$ is a monad. It distributes 
over $\free_-\colon\poly\to\poly$ as a pointed endofunctor but not as a monad.%
\footnote{In particular, the following diagram does not commute when applied to $\yon$:
\[
\begin{tikzcd}[ampersand replacement=\&]
  \freecol\circ\free\circ\free\ar[r]\ar[d]\&
  \free\circ\freecol\circ\free\ar[r]\&
  \free\circ\free\circ\freecol\ar[d]\\
  \freecol\circ\free\ar[rr]\&\&
  \free\circ\freecol\ar[ull, phantom, "\text{Does not commute}"]
\end{tikzcd}
\]
}

Suppose that $m,n$ are $(\otimes)$-monoids, that $c$ is a $(\otimes)$-comonoid, and that $p$ is any polynomial; then:
\begin{align}
	m\otimes n&\quad\tn{carries the structure of a $(\otimes)$-monoid}\\
	m\times n&\quad\tn{carries the structure of a $(\otimes)$-monoid}\\
	m\tri n&\quad\tn{carries the structure of a $(\otimes)$-monoid}\\
	[p, m\otimes p]&\quad\tn{carries the structure of a $(\otimes)$-monoid}\\
	[c\otimes p, p]&\quad\tn{carries the structure of a $(\otimes)$-monoid}
\end{align}

The $(\otimes)$-comonoids $(p,\epsilon,\delta)$ are ``sets of monoids'', i.e.\ for each position $I:p(1)$ the direction set $p[I]$ carries a monoid structure, where $\epsilon$ provides the unit and $\delta$ provides the multiplication.

If $(p,\epsilon,\delta)$ is a $\otimes$-comonoid, the functor $p\otimes-$ induces a colax monoidal functor
\begin{equation}
	(p\otimes -)\colon(\poly,0,\vee)\to(\poly,0,\vee)
\end{equation}
given by $(p\otimes 0)\cong0$ and $p\otimes (q_1\vee q_2)\cong(p\otimes q_1)+(p\otimes q_1\otimes q_2)+(p\otimes q_2)\to(p\otimes q_1)\vee(p\otimes q_2)$.

If $(m,\eta,\mu)$ is a $(\times)$-monoid, then for any $p:\poly$, there is a $\otimes$-comonoid structure
\begin{equation}
	\epsilon\colon\lens{p}{m}\to\yon
	\qqand
	\delta\colon\lens{p}{m}\to\lens{p}{m}\otimes\lens{p}{m}
\end{equation}
on $\lens{p}{m}$ given by the transposes of the following maps
\begin{gather}
	p\to \yon\to m=\yon\tri m\\
	p\to pp\to\left(\lens{p}{m}\tri m\right)\left(\lens{p}{m}\tri m\right)\to\left(\lens{p}{m}\otimes\lens{p}{m}\right)\tri mm\to\left(\lens{p}{m}\otimes\lens{p}{m}\right)\tri m
\end{gather}

For any $(\times)$-monoid $(m,\eta,\mu)$, \eqref{eqn.lax_times_tensor} induces a $(\otimes)$-monoid structure on $m\yon$
\begin{equation}
	\yon\cong1\yon\To{\eta\yon} m\yon
	\qqand
	m\yon\otimes m\yon\to (mm)\yon\To{\mu\yon} m\yon.
\end{equation}

\paragraph{$(\tri)$-monoids and comonoids.} These will be discussed in \cref{chap.monad_comonad}.

\chapter{Polynomial monads and comonads}\label{chap.monad_comonad}

A polynomial monad is a monoid in $(\poly,\yon,\tri)$, i.e.\ a tuple $(t,\eta,\mu)$ where $t:\poly$ and where $\eta\colon\yon\to t$ and $\mu\colon t\tri t\to t$ satisfy the monoid laws. Note that the maps $\eta$ and $\mu$ are often asked to be cartesian; we refer to these as \emph{cartesian monads}. Finitary cartesian monads can be identified with $\Sigma$-free (one-object) operads: $t(1)$ is the set of operations and for each $T:t(1)$ the set $t[T]$ is its arity; the unit and multiplication $(\eta,\mu)$ correspond to the identity and composition respectively. Morphisms between $(\tri)$-monoids are functors between operads.

A polynomial comonad is a comonoid in $(\poly,\yon,\tri)$, i.e.\ a tuple $(c,\epsilon,\delta)$ where $c:\poly$ and where $\epsilon\colon c\to\yon$ and $\delta\colon c\to c\tri c$ satisfy the comonoid laws. Polynomial comonads can be identified with categories: $c(1)$ is the set of objects, and for each $C:c(1)$ the set $c[C]$ is the set of outgoing arrows from that object; the counit corresponds to the identities and the comultiplication corresponds to the codomain and composition information. Morphisms of $(\tri)$-comonoids are not functors, but retrofunctors: we leave this to the reader to investigate.

The category of $(\tri)$-comonoids has finite sums and products, with units $0$ and $\yon$ respectively; it also has another monoidal product $\otimes$ with unit $\yon$, which agrees with the usual cartesian product of categories. The forgetful functor from $\tri$-comonoids to $\poly$ is strong monoidal with respect to $(0,+)$ and $(\yon,\otimes)$; we will discuss its right adjoint $\cofree$ below.

The category of $(\tri)$-monoids has products; the terminal object is $(1,!,!)$, and the product $t\times t'$ has unit $(\eta,\eta')$, and its multiplication arises from \eqref{eqn.comp_times} and projection. The forgetful functor from $\tri$-monoids to $\poly$ is strong monoidal with respect to $\times$. 

For any polynomial monad $(t,\eta,\mu)$, and polynomial $p\colon$, the left Kan extension $\lens{p}{p\tri t}$ is a $\tri$-comonad
\begin{equation}\label{eqn.comonad_from_monad}
  \lens{p}{p\tri t}\to\yon
  \qqand
  \lens{p}{p\tri t}\to\lens{p}{p\tri t}\tri\lens{p}{p\tri t}
\end{equation}
As a category, it is equivalent to the full subcategory of the opposite of the Kleisli category $\smset_t$ spanned by the sets $p[I]$ for $I:p(1)$. For example, the Lawvere theory corresponding to $t$ is the $\tri$-comonad (category)
\begin{equation}
	\text{Law}(t)\cong\lens{\uu}{\uu\tri t}
\end{equation}
where $\uu$ is as in \eqref{eqn.list}. The construction in \eqref{eqn.comonad_from_monad} comes from an oplax monoidal functor
\begin{equation}
	\lens{p}{p\tri\,-}\colon(\poly\op,\yon,\tri)\too(\poly,\yon,\irt),
\end{equation}
whose unitor and co-productor have the form
\begin{equation}
	\lens{p}{p\tri\yon}\to\yon
  \qqand
	\lens{p}{p\tri q_1\tri q_2}\to\lens{p}{p\tri q_2}\tri\lens{p}{p\tri q_1}
\end{equation}

For any cardinality $\kappa$, let $\kappa\poly$ denote the full subcategory of $\kappa$-small polynomials, and let $u_\kappa\coloneqq\sum_{|A|<\kappa}\yon^A$. For example, if $\kappa=\aleph_0$ is the smallest infinite cardinal then $u_\kappa=\uu$. Note that $u_\kappa$ has the structure of a cartesian monad. Consider the Kleisli category $\kappa\polycart_{-\tri u_\kappa}$, whose objects are $\kappa$-small polynomials and whose maps are cartesian natural transformations $p\To{\tn{cart}} q\tri u_\kappa$. Then there is an isomorphism of categories
\begin{equation}
	\kappa\polycart_{-\tri u_\kappa}\cong\kappa\Cat{Dir}
\end{equation}
where $\Cat{Dir}$ is the arrow category of $\smset$ and $\kappa\Cat{Dir}$ is the full subcategory spanned by maps with $\kappa$-small fibers. In other words, a Dirichlet map (or bundle map) from $p$ to $q$ can be identified with a cartesian map of polynomials $p\to q\tri u_\kappa$.

For any $p:\poly$ and comonad $(c,\epsilon,\delta)$, the polynomial $\lens{p\tri c}{p}$ also carries a comonad structure
\begin{equation}\label{eqn.selection_caty}
  \lens{p\tri c}{p}\to\yon
  \qqand
  \lens{p\tri c}{p}\to\lens{p\tri c}{p}\tri\lens{p\tri c}{p}
\end{equation}
These were called \emph{selection categories} and discussed \href{https://topos.site/blog/2021/12/creating-new-categories-from-old-selection-categories/}{here}. 
For any comonad $c$, this construction constitutes a normal oplax monoidal functor
\begin{equation}
\lens{-\tri c}{-}\colon(\poly,1,\times)\to(\catsharp,\yon,\otimes),
\end{equation}
i.e.\ there is a counitor isomorphism and a coproductor (bijective on objects) retrofunctor
\begin{equation}
 \lens{1\tri c}{1}\cong\yon
 \qqand
 \lens{pq\tri c}{pq}\to\lens{p\tri c}{p}\otimes\lens{q\tri c}{q}.
\end{equation}
The construction in \eqref{eqn.selection_caty} comes from an oplax monoidal functor
\begin{equation}
	\lens{p\tri\,-}{p}\colon(\poly,\yon,\tri)\to(\poly,\yon,\tri)
\end{equation}
whose unitor and co-productor have the form
\begin{equation}
  \lens{p\tri\yon}{p}\to\yon
  \qqand
  \lens{p\tri q_1\tri q_2}{p}\to\lens{p\tri q_1}{p}\tri\lens{p\tri q_2}{p}.
\end{equation}

For any polynomial $p:\polycart$, the functor $(p\upup -)$ from \eqref{eqn.upup_formula} sends comonads to comonads. If $c$ is a comonad with corresponding category $\cat{C}$, then $p\upup c$ is the comonad with corresponding category $\sum_{I:p(1)}\cat{C}^{p[I]}$.

Using the maps \eqref{eqn.dayx} and \eqref{eqn.day+}, we have interesting consequences for exponentiating with base $t$, where $t$ is a monad. First, it gives a surprising lax monoidal functor $(\poly,\yon,\tri)\op\to(\poly,1,\times)$:
\begin{equation}\label{eqn.monad_times}
	1\to t^\yon
	\qqand
	t^p\times t^q\to t^{p\tri q}.
\end{equation}
It also gives a lax monoidal functor $(\poly,1,\times)\op\to(\poly,\yon,\tri)$:
\begin{equation}\label{eqn.monad_tri}
	\yon\to t^1
	\qqand
	t^p\tri t^q\to t^{p\times q}.
\end{equation}

Thus we find a $(\poly,1,\times)$-enrichment of the Kleisli category $\smset_t$ for any polynomial monad $t$. Indeed, between any two sets $A,B:\smset$ the hom-polynomial is $t^{A\yon^B}\cong \yon^A\tri t\tri(\yon+B)$. The underlying category has its homsets given by maps from the monoidal unit, $\poly(1,t^{A\yon^B})\cong\yon^A\tri t\tri B=\smset(A,t(B))$ as desired. The identity and composition maps
\begin{equation}
	1\to t^{A\yon^A}
	\qqand
	t^{A\yon^B}\times t^{B\yon^C}\to t^{A\yon^C}
\end{equation}
arise from \eqref{eqn.monad_times} via the obvious maps $A\yon^A\to\yon$ and $A\yon^C\to A\yon^B\tri B\yon^C$. This embeds fully faithfully into the Kleisli category $\poly_{t\tri\,-}$, which is also enriched in $(\poly,1,\times)$ with hom-polynomials $t^{\tiny\lens{p}{q}}$.

From \eqref{eqn.monad_tri} it follows that $t^p$ is a monad for any $p:\poly$. 
\begin{equation}\label{eqn.t^p_monad}
	\yon\to t^p\qqand t^p\tri t^p\to t^p.
\end{equation}
In addition, when $c$ is a comonoid, the polynomial $t^c$ has a monoid structure for several monoidal structures:
\begin{align}
  0\to t^c &\qqand t^c + t^c\to t^c,\\
  \label{eqn.t^c_times_monoid}
  1\to t^c &\qqand t^c\times t^c\to t^c,\\
  \yon\to t^c &\qqand t^c\otimes t^c\to t^c,\\\label{eqn.monad_comonad_power}
  \yon\to t^c &\qqand t^c\tri t^c\to t^c.
\end{align}
In fact, for any polynomial comonad $c$, the functor $-^c$ is a monad on polynomial monads: for any polynomial monad $t$ the maps
\begin{equation}
	t\to t^c
	\qqand
	(t^c)^c\to t^c
\end{equation}
induced by $c\to 1$ and $c\to c\times c$ are in fact monad maps.

For any monad $t$ and comonad $c$, their Dirichlet hom $[c,t]$ also has a monad structure by duoidality:%
\footnote{The isomorphism $[c,t]\tri 1\cong t^c\tri 0$ induces a natural map of polynomials $[c,t]\to t^c$, but it does not respect monad multiplication.}
\begin{align}
	\yon\to[c,t]
	&\qqand
	[c,t]\tri[c,t]\to[c,t].
\end{align}
By this and \eqref{eqn.innerhom_adj}, the polynomial
\begin{equation}
	[S\yon,S\yon]\cong [S\yon^S,\yon]\cong (S\yon)^S
\end{equation}
is a monad, the \emph{$S$-state monad}, for any set $S:\smset$.

For monads, but more generally for pointed endofunctors $\yon\to t$, the map $t\mapsto t^p$ is lax monoidal for any $p:\poly$:
\begin{equation}
	\yon\to\yon^p
	\qqand
	t^p\tri u^p\to(t\tri u)^p.
\end{equation}

Important monads include $\List$, see \eqref{eqn.list}, and $\lott$:
\begin{equation}
	\lott\coloneqq\sum_{N:\nn}\sum_{P:\Delta_N}\yon^N
	\where
	\Delta_N\coloneqq\{P\colon\{1,\ldots,N\}\to[0,1]\mid1=P(1)+\cdots+P(N)\}\}.
\end{equation}
We can build new monads from old as follows. Suppose that $t,t'$ are monads, $c$ is a comonad, $p:\poly$ is any polynomial, and $M:\smset$ is a set. Then
\begin{align}
	tt'&\quad\tn{carries the structure of a monad}\\
	\label{eqn.t^p}
	t^p&\quad\tn{carries the structure of a monad}\\
	\label{eqn.[c,t]}
	[c,t]&\quad\tn{carries the structure of a monad}\\
	t\tri M\yon&\quad\tn{carries the structure of a monad}.
\end{align}
As an example of \eqref{eqn.t^p}, $t^{R\yon^E}\cong\yon^R\tri t\tri(\yon+E)$ is the ``$R$-reader'' monad, followed by $t$, followed by the ``$E$-exceptions'' monad; thus we can think of \eqref{eqn.t^p} as allowing for exception types that depend on what is read in. As an example of \eqref{eqn.[c,t]}, $[S\yon^S,t]\cong\yon^S\tri t\tri S\yon$ is sometimes called the ``$S$-parser'' monad.

There is a cofree $\tri$-comonoid (often called the cofree comonad) construction on $\poly$:
\begin{equation}\label{eqn.cofree}
\begin{tikzcd}[column sep=60pt]
	\Cat{Comon}(\poly)
  	\ar[from=r, shift left=8pt, "\cofree"]
		\ar[from=r, phantom, "\scriptstyle\bot"]
  	&
	\poly
		\ar[from=l, shift right=-8pt, "U"]
\end{tikzcd}
\end{equation}
where $U$ is the forgetful functor that sends a comonoid to its carrier. The cofree comonoid $\cofree_p$ on $p:\poly$ is carried by the limit
\begin{equation}
\cofree_p\coloneqq\lim(\cdots\to p_{n+1}\To{f_n} p_n\to\cdots\to p_1\To{f_0} p_0)
\end{equation}
where the $p_k$ are defined inductively as follows:
\begin{align}
	p_0&\coloneqq\yon&p_{k+1}&\coloneqq (p\tri p_k)\times\yon\\
\intertext{and the maps $f_k\colon p_{k+1}\to p_k$ are defined inductively as follows:}
	p_1=p\times\yon&\To{f_0\coloneqq\tn{proj}}\yon=p_0&p_{k+2}=(p\tri p_{k+1})\times\yon&\To{f_{k+1}\coloneqq(p\tri f_k)\times\yon}(p\tri p_{k})\times\yon=p_{k+1}
\end{align}
The map $\cofree\to\yon$ is easy and the map $\cofree\to\cofree\tri\cofree$ is given by maps $p_{m+n}\to p_m\tri p_n$, which themselves arise by induction on $n$, properties of $\cocl{}$, and maps $p\tri p_m\to p_m\tri p$ that  arise by induction on $m$. There is an isomorphism of polynomials
\begin{equation}\label{eqn.cofree_iso}
	\cofree_p\To{\cong} (p\tri\cofree_p)\times\yon.
\end{equation}
If $p\to q$ is cartesian, so is $\cofree_p\to\cofree_q$. In particular, \eqref{eqn.cofree_iso} gives rise to Lambek's lemma, which says that the coalgebra map for the finaal algebra $\cofree_p\tri 1$ is a bijection
\begin{equation}
		\cofree_p\tri 1\To\cong p\tri(\cofree_p\tri 1).
\end{equation}

In many different ways, the cofree comonad functor $\cofree\colon\poly\to\poly$ is lax monoidal as it maps out of the $(\yon,\otimes)$ monoidal structure:%
\footnote{Recall from \eqref{eqn.ovee} that $p\ovee q\coloneqq \sum_{(I,J): p(1)\times q(1)}\yon^{p[I]+p[I]\times q[J]+q[J]}$.}
\begin{align}
  \cofree_p\otimes\cofree_q&\to\cofree_{p\times q}\\
  \label{cofree_lax_monoidal}
  \cofree_p\otimes\cofree_q&\to\cofree_{p\otimes q}\\
  \cofree_p\otimes\cofree_q&\to\cofree_{p\tri q}\\
	\cofree_p\otimes\cofree_q&\to\cofree_{p\ovee q}  
\end{align}
It also has natural comonoid homomorphisms of the form
\begin{equation}
	\cofree_{[p,q]}\otimes\cofree_{[p',q']}\to\cofree_{[p+p',q+q']}
\end{equation}
that arise from the counits of any comonoid, as well as the distributivity of $\otimes$ over $+$.


For any $p:\poly$ and $A:\smset$, there is an isomorphism
\begin{equation}
	\cofree_{Ap}
	\cong
	\cofree_p\tri A\yon.
\end{equation}
In other words, $\cofree_{\cofree_A\tri p}\cong\cofree_p\tri\cofree_A$. For example, $\cofree_{A\yon^B}=\yon^{\List(B)}\tri(A\yon)$.

For any cartesian map $p\tocart q$, the induced map $\cofree_p\tocart\cofree_q$ is also cartesian.


There is a free $\tri$-monoid (often called the free monad) construction on $\poly$
\begin{equation}\label{eqn.free_monad}
\begin{tikzcd}[column sep=60pt]
	\poly
  	\ar[from=r, shift left=8pt, "U"]
		\ar[from=r, phantom, "\scriptstyle\bot"]
  	&
	\Cat{Mon}(\poly)
		\ar[from=l, shift right=-8pt, "\free"]
\end{tikzcd}
\end{equation}
where $U$ is the forgetful functor that sends a monoid to its carrier. 
The free $\tri$-monoid on $q$ can be constructed as the colimit:
\begin{equation}
	\free_q\coloneqq\colim_{\alpha<\kappa}(\cdots\from q_{\alpha+1}\From{g_\alpha}q_\alpha\from\cdots\from q_1\From{g_0} q_0)
\end{equation}
where $\kappa$ is any upper bound on the cardinalities of $(q[J])_{J\in q(1)}$ and where the $q_\alpha$ are defined inductively as follows:
\begin{equation}
	q_0\coloneqq \yon
	\qqand
	q_{\alpha+1}\coloneqq \yon+(q\tri q_k)
	\qqand
	q_\beta\coloneqq\colim_{\alpha<\beta}q_\alpha
\end{equation}
where $\beta$ is any limit ordinal, and where the maps $g_k\colon q_{\alpha}\to q_{\alpha+1}$ are defined inductively as follows:
\begin{equation}
	q_0=\yon
	\To{g_0\coloneqq\tn{incl}}\yon+q=q_1
	\qqand
	q_{\alpha+1}=\yon+(q\tri q_{\alpha})
	\To{g_{\alpha+1}\coloneqq\yon+(q\tri g_\alpha)}\yon+(q\tri q_{\alpha+1})=q_{\alpha+2}
\end{equation}
The maps into $q_\beta$ are just the colimit inclusions when $\beta$ is a limit ordinal. Note that each map $g_\alpha$ is cartesian. 

Analogous to \cref{eqn.cofree_iso} there is an isomorphism of polynomials
\begin{equation}
	 \yon+(p\tri\free_p)\To{\cong}\free_p.
	\end{equation}
In particular, this gives rise to Lambek's lemma that the algebra map for the initial algebra $\free_p\tri 0$is a bijection
\begin{equation}
		p\tri(\free_p\tri 0)\To\cong\free_p\tri0.
\end{equation}

For any polynomial $p:\poly$ and set $A:\smset$, there is a natural isomorphism
\begin{equation}
  \free_p\tri(\yon+A)
  \cong
	\free_{p+A}.  
\end{equation}
Said another way, $\free_p\tri\free_A\cong\free_{\free_A\tri p}$. In particular, we can apply $-\tri 0$ to obtain:
\begin{equation}
	\free_p\tri A\cong\free_{p+A}(0)
\end{equation}
and the result $\free_p\tri A$ is isomorphic to the free $p$-algebra on the set $A$.

Using the monad structure \eqref{eqn.monad_tri} on $\free_x^p$, one obtains a tensorial strength on the free monad endofunctor $\free$ on $(\poly,1,\times)$:
\begin{equation}\label{eqn.p_action_free}
	p\times\free_q\to\free_{p\times q}
\end{equation}

The free monad functor $\free\colon\poly\to\poly$ is not lax monoidal with respect to $\otimes$ on both sides,%
\footnote{For example, there is no map of polynomials 
$
  \free_1\otimes\free_0\cong\yon+1
  \To{??}
  \yon\cong\free_0\
$.
}
but in many different ways, the free monad \emph{functor} is lax monoidal with the $(0,\OR)$ monoidal structure from \eqref{eqn.vee}.
\begin{align}
\label{eqn.vee1}
	\free_p+\free_q&\to\free_{p\OR q}\\
\label{eqn.vee3}
	\free_p\tri\free_q&\to\free_{p\OR q}\\
\label{eqn.vee2}
	\free_p\otimes\free_q&\to\free_{p\OR q}\\
\label{eqn.vee4}
	\free_p\OR\free_q&\to\free_{p\OR q}
\end{align}
Each of the maps \eqref{eqn.vee1} to \eqref{eqn.vee4} is cartesian.

\begin{warning}\label{warn.noncommutative}
Note that the monad multiplication $\free_{\free_-}\to\free_-$ is \emph{not} a monoidal natural transformation with respect to $\OR$: the diagram
\begin{equation}
\begin{tikzcd}[column sep=40pt, row sep=23pt]
	\free_{\free_p}\OR\free_{\free_q}\ar[r, "\varphi_{\free_p,\free_q}"]\ar[d,"\mu_p\OR\mu_q"']&
	\free_{\free_p\OR\free_q}\ar[r, "\free_{\varphi_{p,q}}"]&
	\free_{\free_{p\OR q}}\ar[d, "\mu_{p\OR q}"]\\
	\free_p\OR\free_q\ar[rr, "\varphi_{p,q}"']&\ar[u, phantom, pos=.6, "\text{does not commute}"]&
	\free_{p\OR q}
\end{tikzcd}
\end{equation}
\end{warning}

There are (cartesian) lax monoidal structures factoring \cref{eqn.vee1,eqn.vee3}:
\begin{align}
	\free_p+\free_q&\to\free_{p+q}\\\label{eqn.free_tri_plus}
	\free_p\tri\free_q&\to\free_{p+q}.
\end{align}
The latter induces a (cartesian) lax monoidal structure
\begin{equation}
	\free_p\otimes\free_q\to\free_{p+q}\label{eqn.free_plus_otimes}
\end{equation}
coming from $\free_p\otimes\free_q\To{\indep} \free_p\tri \free_q\To{\eqref{eqn.free_tri_plus}} \free_{p+q}$. However, similarly to \cref{warn.noncommutative}, be warned that it does not play well with duoidality \eqref{duoidal}: the two maps $(\free_p\tri\free_p)\otimes(\free_q\tri\free_q)\tto\free_{p+q}$ do not agree.


Using $\garner{\ovee}$ from \cref{eqn.garner_ovee}, there is also a lax monoidal structure
\begin{equation}
	\free_p\garner{\ovee}\free_q\to\free_{p\garner{\ovee}q}.
\end{equation}

For any monad $t$ and set $A$, there is a distributive law of the form
\begin{equation}
	\free_A\tri t\to t\tri\free_A.
\end{equation}
If $t$ is (i.e.\ its unit and multiplication maps are) cartesian, then there is also a distributive law of the form
\begin{equation}
	t\tri\free_A\to \free_A\tri t,
\end{equation}
e.g.\ $t+1$ is also a monad.

For any $p:\poly$ there is a function
\begin{equation}
	\free_p\tri 0\to\cofree_p\tri 1
\end{equation}
since the initial $p$-algebra $\free_p\tri 0$ is also a $p$-coalgebra and $\cofree_p\tri 1$ is the final $p$-coalgebra. In case $p=A\yon+B$, this extends to a bijection 
\begin{equation}
	\cofree_{A\yon+B}\tri 1
	\cong
	\left(\free_{A\yon+B}\tri 0\right) +\left(\cofree_{A\yon}\tri 1\right)
	\cong\sum_{N:\nn}A^NB+A^\nn.
\end{equation}

For any $p:\polycart$ there is a natural map
\begin{equation}
	(\cofree_p)_*\to(\free_p)_*
\end{equation}
making the following diagram commute:
\begin{equation}
\begin{tikzcd}
	(\cofree_p)_*\ar[r]\ar[d, "\epsilon"']&(\free_p)_*\ar[d, "\epsilon"]\\
	\yon\ar[r, equal]&\yon
\end{tikzcd}
\end{equation}
where for an arbitrary $q:\poly$ we write $\epsilon\colon q_*\to\yon$ to denote the counit map of \cref{eqn.pointing_comonad}.


For any polynomials $p,q:\poly$ and map $\varphi\colon q\tri p\to p\tri q$, there is an induced map
\begin{equation}\label{eqn.spooling}
	\cofree_{p\tri q}\to\cofree_p\tri\free_q
\end{equation}
extending to a comonad map $\lens{\cofree_{p\tri q}}{\free_q}\to\cofree_p$.%
\footnote{The map \eqref{eqn.spooling} was discussed in this \href{https://topos.institute/blog/2023-04-24-spooling-syntax-from-behavior/}{blog post}. To obtain the comonad map $\lens{\cofree_{p\tri q}}{\free_q}\to\cofree_p$, it suffices to check that $c\coloneqq\lens{\cofree_{p\tri q}}{\free_q}$ is indeed a comonad, because there is an obvious map $\cofree_{p\tri q}\to p\tri q\to p\tri\free_q$. To see that $c$ is a comonad, note that $\varphi$ induces a map $q\tri p\tri q\to p\tri q\tri q$ and see \cite[Example 5.11]{spivak2025categories}.}

\chapter{Adjunctions, monads, and comonads on $\poly$}\label{chap.adj_mon_com}

There are adjunctions between $\poly$ and $\smset$ and between $\poly$ and $\smset\op$, each labeled by where they send $p:\poly$ and $A:\smset$:
\begin{equation}\label{eqn.adjunctions}
\begin{tikzcd}[column sep=60pt]
  \poly
  	\ar[from=r, shift left=8pt, "A" description]
		\ar[from=r, shift left=-24pt, "A\yon"']&
  \smset
  	\ar[from=l, shift right=24pt, "p(0)"']
  	\ar[from=l, shift right=-8pt, "p(1)" description]
	\ar[from=l, phantom, "\scriptstyle\bot"]
	\ar[from=l, phantom, shift left=16pt, "\scriptstyle\bot"]
	\ar[from=l, phantom, shift right=16pt, "\scriptstyle\bot"]
\end{tikzcd}
\hspace{1in}
\begin{tikzcd}[column sep=60pt]
	\poly
  	\ar[from=r, shift left=8pt, "\yon^A"]
		\ar[from=r, phantom, "\scriptstyle\bot"]
  	&
	\smset\op
		\ar[from=l, shift right=-8pt, "\Gamma(p)"]
\end{tikzcd}
\end{equation}
We write $A$ to denote $A\yon^0$. All the leftward maps in \eqref{eqn.adjunctions} are fully faithful, and all the rightward maps are essentially surjective. The leftward maps from $\smset$ are also rig monoidal (i.e.\ strong monoidal with respect to $+$ and $\otimes$):
\begin{align}
	A\yon+B\yon&\cong(A+B)\yon&
	A\yon\otimes B\yon&\cong(A\times B)\yon\\
	A\yon^0+B\yon^0&\cong(A+B)\yon^0&
	A\yon^0\otimes B\yon^0&\cong(A\times B)\yon^0
\end{align}
The rightward maps to $\smset$ are also distributive monoidal; indeed by \cref{eqn.comp_plus,eqn.comp_times}, the following hold for any $A:\smset$, in particular for $A:\{0,1\}$.
\begin{equation}
	p(A)+q(A)\cong(p+q)(A)
	\qqand
	p(A)\times q(A)\cong(p\times q)(A)
\end{equation}
We denote the idempotent comonad $p\mapsto p(1)\yon$ by $\lin{p}$. It is strong monoidal with respect to $+,\otimes$ and lax monoidal with respect to $\tri$. 

The functor $\Gamma$ can be derived from the cartesian closure \eqref{eqn.cart_cl}
\begin{equation}
  \Gamma(p)=\yon^p\tri 0=\prod_{I:p(1)}p[I]
\end{equation}
It preserves coproducts, since coproducts in $\smset\op$ are products in $\smset$:
\begin{equation}
	\Gamma(p+q)\cong\Gamma(p)\times\Gamma(q)
\end{equation}

We can say more about $\Gamma$ if we package it with $p\mapsto p(1)$; i.e.\ there is an adjunction
\begin{equation}\label{eqn.rectangle}
\begin{tikzcd}[column sep=60pt]
	\poly
  	\ar[from=r, shift left=8pt, "A\yon^B"]
		\ar[from=r, phantom, "\scriptstyle\bot"]
  	&
	\smset\times\smset\op
		\ar[from=l, shift right=-8pt, "{\big(p(1)\,,\,\Gamma(p)\big)}"]
\end{tikzcd}
\end{equation}
The left adjoint above is comonadic. It is also strong monoidal with respect to coproduct and $\otimes$. To say so requires us to mention that $\smset\times\smset\op$ has a coproduct structure and to specify a $\otimes$-structure on $\smset\times\smset\op$; they are given as follows:
\begin{align}
 (A_1,B_1)+(A_2,B_2)&\coloneqq(A_1+A_2\,,\,B_1\times B_2)\\
  (A_1,B_1)\otimes(A_2,B_2)&\coloneqq(A_1\times A_2\,,\,B_1^{A_2}\times B_2^{A_1})
\end{align}
Returning to our point, the left adjoint in \eqref{eqn.rectangle} is (strong) rig monoidal (preserves $+$ and $\otimes$):
\begin{align}
	(p(1),\Gamma(p))+(q(1),\Gamma(q))&\cong((p+q)(1),\Gamma(p+q))\\
	(p(1),\Gamma(p))\otimes(q(1),\Gamma(q))&\cong((p\otimes q)(1),\Gamma(p\otimes q))
\end{align}
It follows that the unit $p\mapsto p(1)\yon^{\Gamma(p)}$ is normal lax monoidal, i.e.\ the unitor $\yon(1)\yon^{\Gamma(\yon)}\cong\yon$ is an isomorphism and for any $p,q:\poly$, there is a map%
\footnote{
Capucci refers to the vertical map \eqref{eqn.nashator} as the ``Nashator'', as it shows up in Nash equilibria. Indeed, given ``moves'' $I:p(1)$ and $J:q(1)$ by two players, and given a ``payoff matrix'' $\gamma\colon \Gamma(p\otimes q)$, this map returns the payoff ``row'' for $p$ given $q$'s move and the payoff ``column'' for $q$ given $p$'s move.
}
\begin{equation}\label{eqn.nashator}
	p(1)\yon^{\Gamma(p)}\otimes q(1)\yon^{\Gamma(q)}\to(p\otimes q)(1)\yon^{\Gamma(p\otimes q)}
\end{equation}

The map $p\mapsto p_*$ from \cref{eqn.p_star} forms an idempotent comonad on $\polycart$. It is strong monoidal with respect to $\otimes,+$. There is an isomorphism $\lin{p}_*\cong \lin{p}$.

The endofunctor $\free_-\colon\poly\to\poly$ from \cref{eqn.free_monad} is itself a monad
\begin{equation}
	p\to\free_p
	\qqand
	\free_{\free_p}\to\free_p
\end{equation}
and the functor $\cofree_-\colon\poly\to\poly$ from \cref{eqn.cofree} is itself a comonad
\begin{equation}
	\cofree_p\to p
	\qqand
	\cofree_p\to\cofree_{\cofree_p}. 
\end{equation}
Moreover, the former is a module over the latter---see \cite{libkind2024pattern}---i.e.\ for any $p,q$ there is a natural map
\begin{equation}\label{module_easy}
  \free_p\otimes\cofree_q\to\free_{p\otimes q}
\end{equation}
satisfying the action laws for the maps from \cref{cofree_lax_monoidal}.%
\footnote{In fact, there is a natural map $\free_{p\vee p'}\otimes\cofree_{q}\to\free_{(p\otimes q)\vee p'}$.}
  It also satisfies the following coherence for $\free$ as a monad and $\cofree$ as a comonad:
\begin{equation}\label{eqn.monad_comonad_coherence}
\begin{tikzcd}
  p\otimes\cofree_q\ar[r]\ar[d]&p\otimes q\ar[d]\\
  \free_p\otimes\cofree_q\ar[r]&\free_{p\otimes q}
\end{tikzcd}
\hspace{.6in}
\begin{tikzcd}
	\free_{\free_p}\otimes\cofree_q\ar[d]\ar[r]&
	\free_{\free_p}\otimes\cofree_{\cofree_q}\ar[r]&
	\free_{\free_p\otimes\cofree_q}\ar[r]&
	\free_{\free_{p\otimes q}}\ar[d]\\
	\free_p\otimes\cofree_q\ar[rrr]&&&
	\free_{p\otimes q}
\end{tikzcd}
\end{equation}
%
%

For any polynomial $p$, the free monad $\free_p$ is a left $\cofree_{\yon+p}$-comodule, i.e.\ there is a left coaction:
\begin{equation}
	\free_p\to\cofree_{\yon+p}\tri\free_p
\end{equation}
which is induced by the composite $\free_p\To{\cong}\yon+p\tri\free_p\to(\yon+p)\tri\free_p$. 

For any $p,q:\poly$ and map $p\tri q\to q\tri p$ there is an induced distributive law
\begin{equation}
  \free_p\tri\cofree_q\to\cofree_q\tri\free_p.
\end{equation}
A special case comes from the identity on $p\tri p$, which gives a well-known distributive law $\free_p\tri\cofree_p\to\cofree_p\tri\free_p$.

For any polynomial $e:\poly$, exponentiation by $e$ is a monad on $\poly$:
\begin{equation}
	p\to p^e
	\qqand
	(p^e)^e\to p^e
\end{equation}
This monad distributes over the free monad monad
\begin{equation}
	\free_{p^e}\to(\free_p)^e
\end{equation}
so $p\mapsto(\free_p)^e$ is also a monad on $\poly$. Similarly, for any polynomial monad $t$, both the functors $p\mapsto t\tri p$ and $p\mapsto p\tri t$ are monads on $\poly$. The free monad monad distributes over both these monads
\begin{equation}
	t\tri\free_p\to\free_{t\tri p}
	\qqand
	\free_p\tri t\to\free_{p\tri t}
\end{equation}
so $p\mapsto \free_{t\tri p}$ and $p\mapsto \free_{p\tri t}$ are both monads on $\poly$.

For any polynomial comonad $(c,\epsilon,\delta)$, there is a comonad on $\poly$ given by $p\mapsto p\tri c$. By \cref{eqn.comp_plus,eqn.comp_times}, this comonad is strong monoidal with respect to $+$ and $\times$:
\begin{align}
	(p\tri c)+(q\tri c)&\cong (p+q)\tri c&0\tri c&\cong0\\
	(p\tri c)\times(q\tri c)&\cong (p\times q)\tri c&1\tri c&\cong1	
\end{align}
While this is basically trivial, it is useful e.g.\ for $c=2\yon$ and other linear polynomials, because $p\tri 2\yon\cong\sum_{I:p(1)}\sum_{U\ss p[I]}\yon^{p[I]}$ allows one to ``specify a desired outcome as part of the position''.

\chapter{*-Bifibration over $\smset$ and factorization systems}\label{chap.bifib}

The functor
\begin{equation}\label{eqn.bifib}
\big(p\mapsto p(1)\big)\colon\poly\to\smset
\end{equation}
is a *-bifibration. In particular, for any function $f\colon A\to B$, there is an adjoint triple $f_!\dashv f^*\dashv f_*$,
\begin{equation}
\begin{tikzcd}[column sep=50pt]
	\poly_A
		\ar[r, shift left=16pt, "f_!"]
		\ar[r, shift right=16pt, "f_*"']
		\ar[from=r, "f^*" description]
		\ar[r, phantom, shift left=8pt, "\Rightarrow"]
		\ar[r, phantom, shift right=8pt, "\Leftarrow"]
	&
	\poly_B
\end{tikzcd}
\end{equation}
where $\poly_X$ is the category of polynomials with positions $p(1)=X$. The images under the functors $f_!$ and $f_*$ of $p:\poly_A$ are given by
\begin{equation}
	f_!(p)\coloneqq\sum_{b: B}\yon^{\;\prod\limits_{b=fa}p[a]}
	\qqand
	f_*(p)\coloneqq\sum_{b: B}\yon^{\;\sum\limits_{b=fa}p[a]}
\end{equation}
and the image under the functor $f^*$ of $q:\poly_B$ is given by
\begin{equation}
	f^*(q)\coloneqq\sum_{a: A}\yon^{q[fa]}
\end{equation}
For any $p:\poly_A$ and $q:\poly_B$ there are natural maps
\begin{equation}
	p\to f_!(p)
	\qqand
	f^*(q)\to q.
\end{equation}

A morphism $\varphi\colon p\to q$ can be identified with a diagram of the form
\begin{equation}\label{eqn.poly_map}
\begin{tikzcd}
	p(1)\ar[d, "\varphi_1"']\ar[r, "{p[-]}", ""' name=p]&
	\smset\\
	q(1)\ar[ur, bend right, "{q[-]}"', "" name=q]
	\ar[to=p, from=q-|p, Rightarrow, shorten=3pt, "\varphi^\sharp"]
\end{tikzcd}
\end{equation}
The $p\mapsto p(1)$ bifibration \eqref{eqn.bifib} gives us the terms \emph{vertical, cartesian, \text{and} op-cartesian} for a map $\varphi\colon p\to q$ in $\poly$. That is, taking $f\coloneqq\varphi_1$, we have that $\varphi$ is vertical if it is contained in a fiber of the bifibration \eqref{eqn.bifib}, it is cartesian if $p\to f^*(q)$ is an isomorphism, and it is op-cartesian if $f_!(p)\to q$ is an isomorphism. Here are alternative ways to define these notions: $\varphi$ is 
\begin{itemize}
	\item \emph{vertical} if $\varphi_1\colon p(1)\to q(1)$ is an identity in $\smset$,
	\item \emph{cartesian} if $\varphi^\sharp$ is a natural isomorphism, and
	\item \emph{op-cartesian} if the diagram \eqref{eqn.poly_map} is a right Kan extension.
\end{itemize}
More explicitly, $\varphi$ is cartesian iff for each $I: p(1)$, the function $\varphi^\sharp_I\colon p[I]\to q[\varphi_1I]$ is a bijection. It is op-cartesian if for each $J: q(1)$ the map $q[J]\to\prod\limits_{\varphi_1(I)=J}p[I]$ is a bijection.

Let \emph{essentially-vertical} mean that $\varphi_1$ is a bijection, rather than identity. Then there are at least four factorization systems on $\poly$:
\begin{itemize}
	\item (epi, mono),
	\item (essentially-vertical, cartesian),
	\item (op-cartesian, essentially-vertical),
	\item (taut, essentially-vertical epi).
\end{itemize}
A natural transformation $\varphi\colon p\to q$ is \emph{taut} iff for all monomorphisms of sets $f\colon A\to B$, the naturality square for $f$ and $\varphi$ is a pullback. Par\'{e} showed that if $p,q$ are polynomials then $\varphi$ is taut iff for each $I:p(1)$ the map $\varphi_I^\sharp\colon q[\varphi(I)]\to p[I]$ is surjective. These form the left class of a factorization system whose right class are bijective on positions and injective on directions.

We note that $\poly$ is not a regular category because epimorphisms (and hence image factorizations) are not stable under pullback. Indeed, the notationally-indicated bottom map in the pullback diagram
\begin{equation}
\begin{tikzcd}
	2\yon\ar[r, "!"]\ar[d]&
	\yon^{\{1,2\}}\ar[d, "\psi"]\\
	\yon^{\{a,b\}}+\yon^{\{c,d\}}\ar[r, "\varphi"']&
	\yon^{\{1ac, 1bd, 2ad, 2bc\}}\ar[ul, phantom, very near end, "\lrcorner"]
\end{tikzcd}
\end{equation}
is an epimorphism, but its pullback along $\psi$, the unique map $2\yon\to\yon^2$ is not. 

\section*{Acknowledgments}
\thanksAFOSR{FA9550-20-1-0348 and FA9550-23-1-0376}.

\printbibliography 
\end{document}